\documentclass[10pt]{amsart}
\usepackage{psfrag}
\usepackage{graphicx}

\newtheorem{theorem}{Theorem}[section]
\newtheorem{lemma}[theorem]{Lemma}
\newtheorem{corollary}[theorem]{Corollary}
\newtheorem{proposition}[theorem]{Proposition}

\theoremstyle{definition}

\theoremstyle{remark}
\newtheorem{remark}[theorem]{Remark}

\numberwithin{equation}{section}

\newcommand{\Er}{\mathbb{R}}
\newcommand{\Ce}{\mathbb{C}}

\newcommand{\call}{\mathcal{L}}

\newcommand{\calr}{\mathcal{R}}
\newcommand{\ve}{\varepsilon}

\newcommand{\norm}[2][{}]{\lVert#2\rVert_{{#1}}}
\newcommand{\vnorm}[2][{}]{\left\Vert#2\right\Vert_{#1}}
\newcommand{\abs}[2][{}]{\lvert#2\rvert_{#1}}
\newcommand{\vabs}[2][{}]{\left\vert#2\right\vert_{#1}}

\begin{document}

\title[Corner penalty\ \today]{Penalty approximation for non smooth
constraints in vibroimpact}

\author{L\ae titia Paoli}
\address{MAPLY, UMR 5585 CNRS and Analyse Num\'erique, Saint-Etienne \\
Facult\'e des Sciences,
Universit\'e Jean Monnet \\
23 Rue du Docteur Paul Michelon\\
42023 St-Etienne Cedex 2\\
France}
\email{paoli@anumsun1.univ-st-etienne.fr}

\author{Michelle Schatzman}
\address{MAPLY, UMR 5585 CNRS \\
Universit\'e Lyon 1 \\
69622 Vil\-leur\-ban\-ne Cedex \\
France}
\email{schatz@maply.univ-lyon1.fr}



\begin{abstract}
We examine the penalty approximation of the free motion of a material
point in an angular domain; we choose an over-damped penalty
approximation, and we
prove that if the first impact point is not at the vertex, then, the
limit of the approximation exists and is described by Moreau's rule
for anelastic impacts.

The proofs rely on validated asymptotics and use some classical tools of
the theory of dynamical systems.

\end{abstract}

\maketitle

\section{Introduction}

Mathematical results relative to the convergence of a penalty
approximation of impact problems have been obtained by several authors
when the energy is conserved; see for instance \cite{MR55:3578},
\cite{schatz78}, and also \cite{MR83b:49008}, \cite{MR84i:34016}
\cite{buttazzopercivale}, \cite{MR87m:73025}, \cite{MR90j:70025},
\cite{MR93k:58094a}, \cite{MR93m:70015}, \cite{MR96j:70016},
 \cite{MR95m:70003} and \cite{paoli00}.

When energy may be lost at impact, the convergence of the penalty
approximation has been treated in \cite{pasc93} in the case of a
convex set of constraints with smooth boundary. In this article, we
defined a penalty approximation for which the limit solution satisfies
a Newton condition at impact: the normal component of the velocity is
reversed and multiplied by a restitution coefficient $e\in ]0,1]$ and
the tangential component is transmitted.

So far, we are not aware of any mathematical results on the
convergence of the penalty approximation when the boundary is not
smooth and energy can be lost at impact.

Here, we study the penalty approximation of the motion of a free
particle constrained to stay inside an angular domain of $\Er^2$: we
choose a class of penalty approximations for which the restitution
coefficient vanishes in the limiting problem and we characterize
precisely the limit of the sequence of solutions of the approximated
problem when the first impact does not take place at the corner.

We compare our results to the ones given by the selection rule of
Moreau \cite{moreau1}, and we find complete agreement.

%
%
%
%

\section{The first part of the motion and the mathematical strategy}

Let us describe more precisely the problem and the method of solution.

Given $\bar \theta\in (0, \pi)$, we let $K$ be the set
\begin{equation}
K=\bigl\{(x_1, x_2)\in\Er^2: x_1\le 0 \text{ and } x_1\cos \bar\theta +
x_2\sin\bar\theta \le 0\bigr\}.
\end{equation}

\begin{figure}
\psfrag{erre1}{$\calr_1$}
\psfrag{erre2}{$\calr_2$}
\psfrag{erre3}{$\calr_3$}
\psfrag{bartheta}{$\bar \theta$}
\psfrag{x1}{$x_1$}
\psfrag{x2}{$x_2$}
\psfrag{y1}{$y_1$}
\psfrag{y2}{$y_2$}
\psfrag{Ka}{$K$}
\begin{center}
\includegraphics[height=6 cm]{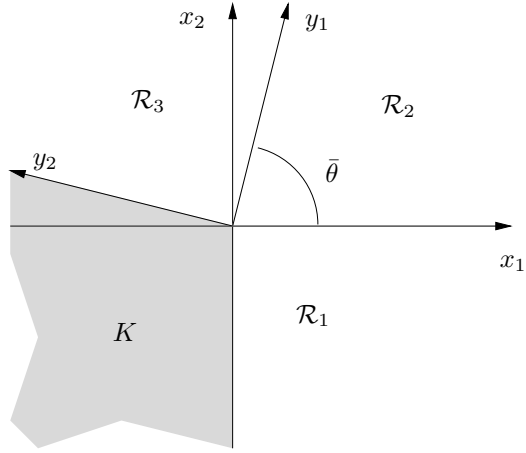}
\bigskip
\caption{The regions and the coordinates in the plane.}\label{fig:region}
\end{center}
\end{figure}

The closure of the complement of $K$ is partitioned into three
regions: 
\begin{align*}
\calr_1&=\{(x_1,x_2): x_1\ge 0 \text{ and } x_2\le 0\},\\
\calr_2&=\{(x_1,x_2): x_2\ge 0 \text{ and } -x_1\sin\bar\theta + x_2
\cos\bar \theta \le 0\}\\
\calr_3&=\{(x_1,x_2): x_1\cos\bar \theta + x_2\sin\bar \theta \ge 0\notag\\
&\qquad\text{ and
} -x_1\sin\bar\theta + x_2 \cos\bar \theta \ge 0\}.
\end{align*}

In each of these regions, the projection onto $K$, which is known to be a
contraction, takes different forms: 
\begin{equation*}
P_K x = \begin{cases}
(0,x_2)^T, & \text{if $x\in \calr_1$,}\\
0, & \text{if $x\in \calr_2$,}\\
(-x_1\sin \bar \theta + x_2\cos\bar \theta)(-\sin\bar \theta, \cos\bar
\theta)^T,& \text{if $x\in \calr_3$.}
\end{cases}
\end{equation*}

The penalty approximation used in \cite{pasc93} is defined as follows:
we define a function $G$ of two arguments $u\in \Er^2$ and $v\in\Er^2$ by
\begin{equation*}
G(u,v)=\begin{cases}
\frac{\displaystyle\bigl(v\cdot (u-P_Ku)\bigr)(u-P_Ku)}{\displaystyle |u-P_Ku|^2},&\text{if $u\notin K$,}\\
0,&\text{otherwise.}
\end{cases}
\end{equation*}

Then, the penalized approximation to the impact problem in $K$, in the
absence of exterior forces is
given by
\begin{equation}
\label{approx}
\ddot u_k + 2 \alpha \sqrt k G(u_k,\dot u_k) + k (u_k-P_Ku_k)=0.
\end{equation}

In this formulation, the number $k$ is the stiffness of the spring
which describes the reaction of the wall, and the choice of the scale
$\sqrt k$ is the exact choice which ensures convergence in the smooth
case, as $k$ tends to infinity. Here, we choose $\alpha >1$: it is the
over-damped choice and it 
will lead to a vanishing restitution coefficient as we shall see. 

%

Take initial conditions given by
\begin{equation}
\label{initial} 
\begin{split}
&u_k(0)=\bigl(0, x_2(0)\bigr)^T, \quad x_2(0)<0, \\
&\dot u_k(0)= \bigl(\dot
x_1(0), \dot x_2(0)\bigr)^T,\quad \dot x_1(0)>0, \dot x_2(0)>0.
\end{split}
\end{equation}

These initial conditions mean that at time $t=0$, the representative
point of the system is on the boundary of $K$, in region $\calr_1$, and
outgoing as well as taking the direction of the corner.
In particular, this choice of initial conditions means that the first
impact time is $t=0$. 

The roots of the characteristic equation 
\begin{equation*}
\xi^2 + 2\alpha \xi + 1=0
\end{equation*}
of the over-damped equation
\begin{equation}
\label{overdamped}
\ddot y + 2 \alpha \dot y + y=0
\end{equation}
are given by the formulas
\begin{equation}
\label{roots}
\Delta= \alpha^2 -1, \quad \xi_1 =-\alpha +\sqrt \Delta, \quad \xi_2 =
-\alpha -\sqrt \Delta.
\end{equation}
Both $\xi_1$ and $\xi_2$ are negative.

As long as the representative point of the system lies in $\calr_1$,
we perform the change of variables
$r(t)=x_1(t)\ge 0$ and
$s(t)=x_2(t)\le 0$. In these new coordinates, 
\eqref{approx} becomes the decoupled system
\begin{align*}
&\ddot r + 2\alpha\sqrt k \dot r + kr=0,\\
&\ddot s=0.
\end{align*}
Its solution is given explicitly by
\begin{align}
r(t)&=\frac{\dot r(0)}{2\sqrt{\Delta k}} \left (e^{\xi_1t\sqrt k}-
e^{\xi_2 t\sqrt k}\right),\label{eqr}\\
s(t)&=s(0)+ t\dot s(0).\label{eqs}
\end{align}
For all positive $t$, $r(t)$ given by \eqref{eqr} remains strictly
positive; $s$ reaches the value $0$ at the time
\begin{equation}
\label{t0}
t_0=-s(0)/\dot s(0).
\end{equation}
Therefore, at the boundary between regions $\calr_1$ and $\calr_2$ we have
\begin{equation}
\begin{split}
&s(t_0-0)=0,\quad \dot s(t_0-0)=\dot s(0), \\ 
& r(t_0-0)= \frac{\dot r(0)}{2\sqrt{\Delta k}}\left (e^{\xi_1t_0\sqrt k}-
e^{\xi_2 t_0\sqrt k}\right), \\
& \dot r(t_0-0)=\frac{\dot r(0)}{2\sqrt{\Delta }}\left (\xi_1e^{\xi_1t_0\sqrt k}-
\xi_2e^{\xi_2 t_0\sqrt k}\right).
\end{split}
\label{init0}
\end{equation}

In order to study the motion in region $\calr_2$, we use polar
coordinates i.e. $u_k=r e^{i\theta}$ and we define scaled
functions and variables $R$, $\Theta$ and $\tau$ by
\begin{equation*}
r(t)= \eta R(\tau)/\sqrt k, \quad \tau =(t-t_0)\sqrt k, \quad \eta=
e^{\xi_1 t_0 \sqrt k\, /2}, \quad \Theta(\tau)=\theta(t).
\end{equation*}

\begin{figure}
\psfrag{derre}{$\dot R$}
\psfrag{erre}{$R$}
\psfrag{errec}{$R_c$}
\psfrag{asympto1}{$A_1$}
\psfrag{asympto2}{$A_2$}
\centering
\includegraphics[width=0.79\textwidth]{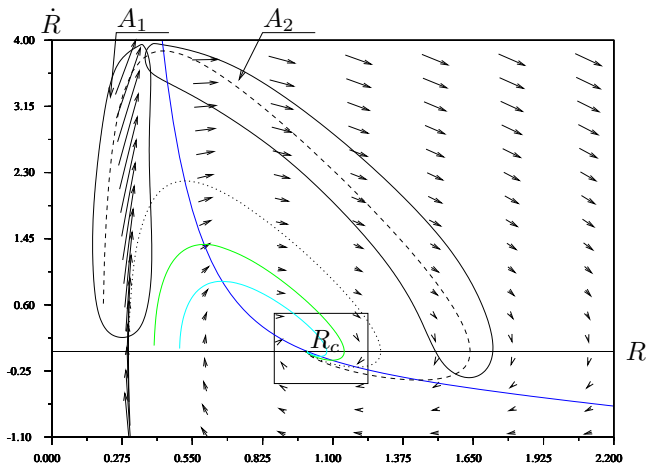}

\caption{The phase portrait for the $R$ equation, and several
trajectories of solutions in the $R, \dot R$ plane. $A_1$: region of
the first asymptotic (section \ref{sec:earliest}); $A_2$: region of the second
asymptotic (section \ref{sec:second}).}\label{fig:phase}
\end{figure}

We have represented in Fig. \ref{fig:phase}
the numerically computed trajectories (dotted or
dashed lines) and the vector field of the ordinary differential
equation for $R$ and $\dot R$.

In the new variables, the system under consideration becomes
\begin{equation}
\label{coin:eq:6}
\ddot R -\frac{E(1-\ve)^2}{R^3} + 2\alpha \dot R +  R=0.
\end{equation}
with 
\begin{equation}
\dot \Theta =\frac{\sqrt E\, (1-\ve)}{R^2},\label{coin:eq:10}
\end{equation}
and the detailed derivation of these equations is performed in
subsection \ref{sec:deriv-scal-equat}.
In equations \eqref{coin:eq:6} and \eqref{coin:eq:10}, $\ve=o(1)$ and
$E$ is a fixed number depending only on the initial conditions and
$\alpha$. The
representation given in Fig.\ref{fig:phase} will help us to explain
how the solution of \eqref{coin:eq:6} behaves, with appropriate
consequences on the angle $\Theta$.

In region $A_1$, $R$ decreases somewhat and then increases, $\dot R$
increases from a size equivalent to $C\eta$ to a size equivalent to
$C/\eta$ in that same region. The dominant terms in equation
\eqref{coin:eq:6} are $\ddot R$ and $E(1-\ve)^2/R^3$; therefore, we
are led to the problem
\begin{equation}
\label{ansatz0}
\ddot R_1 -\frac{E}{R_1^3}=0, \quad R_1(0)=R(0),\quad \dot R_1(0)=\dot R(0).
\end{equation}
We study the solution $R_1$ of \eqref{ansatz0}
in subsection \ref{sec:ansatz-motion-calr_2}, as well as the evolution
of the function  $\Theta_1$ satisfying 
\begin{equation*}
\dot \Theta_1= \sqrt E/ R_1^2;
\end{equation*}
this can be done explicitly, thanks to the simple structure of
\eqref{ansatz0}. In subsection \ref{sec:kern-line-equat}, we study the
kernel of the linearized \eqref{ansatz0} at $R_1$, as a preparation for
the validation of this first asymptotic, a task which is completed in
\ref{sec:valid-earl-asympt} on the interval $[0, \tau_1]$, where 
$\tau_1$ is equal to $\eta^{\gamma_1}$, with $\gamma_1$ belonging to
$(1,2)$. We conclude section \ref{sec:earliest} by 
Proposition \ref{thr:2} which shows that $R$ is equivalent to $R_1$
over $[0, \tau_1]$ and $\dot R$ is equivalent to $\dot R_1$ over
$[\eta^3, \tau_1]$. 
The proof is basically a consequence of the fixed point theorem with a
number of technical estimates.

In Section \ref{sec:case-bar-thetapi2}, assuming $\bar \theta <\pi/2$,
we are able to exploit the above equivalents and to prove that
$\Theta$, solution of \eqref{coin:eq:10}, crosses through $\bar
\theta$ at some time $\bar \tau < O(\eta^2)$. Moreover, our estimates
enable us to describe the limit $u_\infty$ of $u_k$ as $k$ tends to
infinity.
Let $\Pi_1$ be the orthogonal projection on $\{x_1=0\}$, and let
$\Pi_2$ be the orthogonal projection on $\{x_1\cos \bar \theta + x_2
\sin \bar \theta =0\}$; then
\begin{equation*}
u_\infty(t) =
\begin{cases}
u(0) + t\Pi_1 \dot u(0) &\text{ if $0\le t \le t_0$,}\\
(t-t_0)\Pi_2 \Pi_1 \dot u(t_0) &\text{ if $t_0\le t$.}
\end{cases}
\end{equation*}

If $\bar \theta \ge \pi/2$, the representative point of the system
enters region $A_2$ of Fig. \ref{fig:phase}. We have to produce an
asymptotic for the solution of \eqref{coin:eq:6}; in this region, it
is the linear part of this ordinary differential equation which is
dominant; more precisely, let $R_2$ be the solution of
\begin{equation*}
\ddot R_2 + 2\alpha \dot R_2 + R_2 =0,
\end{equation*}
with $R_2$ and $\dot R_2$ respectively coinciding with $R$ and $\dot
R$ at time $\tau_1 = \eta^{\gamma_1}$, where, now the interval of
$\gamma_1$ is reduced to $(1,4/3)$. 

The validation of this ansatz is another consequence of the fixed
point theorem for strict contraction, together with a number of
technical estimates.

Finally, we use classical methods for dynamical systems and prove that
the representative point of the system tends to $(R_c,0)$ as time
tends to infinity: $R_c$ is a number which depends only on the initial
conditions, $\alpha$ and $\ve$. We combine the use of a Lyapunov
functional and some elementary properties of the system to conclude
that $R$ remains bounded from above and away from $0$ for all time
after leaving $A_2$. Observe that the Lyapunov functional gives scant
information in regions $A_2$ and $A_1$: there it takes values of order
$1/\eta$. 

With some technicalities in the case $\bar \theta=\pi/2$, it is
possible to conclude that $\Theta(\tau)$ crosses $\bar \theta$ at some
time $\bar \tau$ and to obtain precise equivalents for $R$, $\dot R$
and $\dot \Theta$ at time $\bar \tau$. After this time, the
representative point of the system \eqref{approx} enters region
$\calr_3$, and we conclude by Theorem \ref{coin:thr:1} that the limit
$u_\infty$ of $u_k$ is given by
\begin{equation*}
u_\infty(t)=\begin{cases}
u(0) + t\Pi_1 \dot u(0) &\text{if $0 \le t \le t_0$,}\\
0 &\text{if $t_0 \le t$.}
\end{cases}
\end{equation*}

Moreau's rule is described as follows: at impact, the outgoing
velocity is projected onto the tangent cone to the convex of
constraints, and the motion proceeds with this new velocity. Thus, it
can be seen that the over-damped penalty approximation agrees
completely with Moreau's rule if the first impact does not take place
at the corner, or very close to it, i.e. at a distance $O\bigl(1/\sqrt
k\bigr)$ from it.

We conjecture that the behavior described here still holds if there is
a right hand side, and the convex is replaced by a set with convex
corners, and smooth and not necessarily convex curves between corners. We
also conjecture that the behavior of the limit of the over-damped penalized
solution is the same in higher spatial dimension.

\section{Equations of the motion around the corner: the earliest
asymptotic}\label{sec:earliest} 

\subsection{Derivation of the scaled equation in
$\calr_2$}\label{sec:deriv-scal-equat} 

After time $t_0$, we arrive into region $\calr_2$, in which it is
convenient to identify $\Er^2$ and $\Ce$ and to use polar coordinates,
i.e.
\begin{equation}\label{polar}
u_k= r e^{i\theta}.
\end{equation}
By continuity, the limits of $r(t)$ from the right and from the left
as $t$ tends to $t_0$ are identical; therefore:
\begin{equation}
\label{inir+}
r(t_0+0)=  \frac{\dot r(0)}{2\sqrt{\Delta k}}\left (e^{\xi_1t_0\sqrt k}-
e^{\xi_2 t_0\sqrt k}\right).
\end{equation}
At $t=t_0$, 
\begin{equation}
\label{initth+}
\theta(t_0+0)=0.
\end{equation}
We differentiate once \eqref{polar} with respect to time:
\begin{equation*}
\dot u_k=\dot r e^{i\theta} + i\dot \theta r e^{i\theta}, 
\end{equation*}
hence
\begin{align}
\dot r(t_0+0) &= \frac{\dot r(0)}{2\sqrt{\Delta }}
\left (\xi_1e^{\xi_1t_0\sqrt k}-
\xi_2e^{\xi_2 t_0\sqrt k}\right),\label{initdr0}
\\
\dot \theta(t_0+0)&=\frac{\dot s(t_0-0)}{r(t_0-0)}
= \frac{\dot s(0) 2\sqrt{\Delta k}}{\dot r(0)\left (e^{\xi_1t_0\sqrt k}-
e^{\xi_2 t_0\sqrt k}\right) }.\label{initdt0}
\end{align}

Let us derive the differential equations satisfied by $r$ and $\theta$.
In region $\calr_2$, the definition of the projection $P_K$ implies that
\eqref{approx} can be written as
\begin{align}
&r\ddot \theta + 2 \dot r \dot \theta =0,\label{eqth}\\
&\ddot r - r\dot \theta^2 + 2\alpha \sqrt k \dot r + k r=0.\label{eqr2}
\end{align}
The motion has central acceleration, therefore it has a first
integral: there exists a constant $\Gamma$ such that
\begin{equation}
\label{invar}
\bigl(r^2\dot\theta\bigr)(t)= \bigl(r^2\dot\theta\bigr)(t_0+0)=\Gamma,
\end{equation}
and, according to \eqref{initdr0} and \eqref{initdt0}, the value of $\Gamma$ is given by
\begin{equation}
\label{valGamma}
\Gamma= r(t_0+0)^2 \frac{\dot s(t_0-0)}{r(t_0+0)}=\frac{\dot r(0)\dot
s(0)}{2\sqrt{\Delta k}} \left (e^{\xi_1t_0\sqrt k}-
e^{\xi_2 t_0\sqrt k}\right).
\end{equation}
We substitute $\dot \theta = \Gamma /r^2$ into \eqref{eqr2}, and we find
the equation in $r$
\begin{equation}
\label{eqrnsc}
\ddot r - \frac{\Gamma^2}{r^3} + 2 \alpha \sqrt k \dot r + kr=0.
\end{equation}
Let us find now appropriate changes of scale which will help us
understand the behavior of $r$ while $u_k$ remains in region
$\calr_2$. An important scale is defined by the number
\begin{equation}
\label{valeta}
\eta=e^{\xi_1t_0\sqrt k \,/ 2},
\end{equation}
which is very small since $\xi_1$ is strictly negative.
We perform the following changes of variables:
\begin{equation}
\label{change}
\tau=(t-t_0)\sqrt k, \quad r(t) = \eta  R(\tau)/\sqrt k.
\end{equation}
In these new variables, we have
\begin{equation*}
\dot r(t)=\eta \dot R(\tau), \quad \ddot r(t)=\eta \sqrt k\ddot R(\tau),
\end{equation*}
so that equation \eqref{eqrnsc} becomes, after a division by
$\eta\sqrt k$
\begin{equation}
\label{eqrnsc1}
\ddot R -\frac{k\Gamma^2}{\eta ^4 R^3} + 2\alpha \dot R +  R=0.
\end{equation}
Let us define
\begin{align}
E&=\frac{\dot r(0)^2 \dot s(0)^2}{4\Delta},\label{valE}\\
\ve&=e^{(\xi_2-\xi_1)t_0\sqrt k}.\label{valeps}
\end{align}
With these notations, we find that $\Gamma$, given by \eqref{valGamma}, is equal to
\begin{equation}
\label{valGamma1}
\Gamma=\frac{\eta^2(1-\ve)\sqrt E}{\sqrt k},
\end{equation}
so that \eqref{eqrnsc1} can be rewritten as
\begin{equation}
\label{eqrnsc2}
\ddot R -\frac{E(1-\ve)^2}{R^3} + 2\alpha \dot R + R=0.
\end{equation}
The initial values for \eqref{eqrnsc2} are given by
\begin{equation}
\label{initR}
R(0)=\frac{\eta\dot r(0)(1-\ve)}{2\sqrt \Delta},
\end{equation}
and
\begin{equation}
\label{initdR}
\dot R(0)=\frac {\eta\dot
r(0)\xi_1(1-\ve\xi_2/\xi_1)}{2\sqrt \Delta}.
\end{equation}

\subsection{Ansatz for the motion in
$\calr_2$}\label{sec:ansatz-motion-calr_2} 

We shall use now an ansatz, namely, we state that the essential part
in the left hand side of \eqref{eqrnsc2} is $\ddot R - E/R^3$.
This comes from the fact that at time $t_0$, $E/R^3$ is very large
with respect to $R$ and $\dot R$ as can be checked from \eqref{initR}
and \eqref{initdR}.
Therefore, we first solve explicitly the equation
\begin{equation}
\label{ansatz}
\ddot R_1 -\frac{E}{R_1^3}=0, \quad R_1(0)=R(0),\quad \dot R_1(0)=\dot R(0).
\end{equation}
If we multiply \eqref{ansatz} by $\dot R_1$ and integrate, we find
that
\begin{equation}
\label{1integ}
\dot R_1^2 + \frac{E}{R_1^2} = \dot R(0)^2 + \frac{E}{R(0)^2}.
\end{equation}
Let us denote by
\begin{equation}
\label{valE1}
W=\dot R(0)^2 + \frac{E}{R(0)^2}
\end{equation}
the value which appears on the right hand side of \eqref{1integ}.
Thanks to the change of variable
$\rho = R_1^2$, equation \eqref{1integ} becomes:
\begin{equation}
\label{eqrho}
\frac{W\dot \rho}{2\sqrt{W\rho -E}} = \pm W.
\end{equation}
At the initial time, $\dot \rho(0)$ is strictly negative, so that in
\eqref{eqrho}
we choose the minus sign on the right hand side, we integrate until a
time $\tau_0$ such that $\dot \rho$ vanishes, and we find
\begin{equation}
\label{solrho-}
\rho(\tau)=\frac{E + (W \tau +\dot R(0) R(0))^2}{W},\quad 0 \le \tau\le
\tau_0.
\end{equation}
The value of $\tau_0$ is given by
\begin{equation}
\tau_0= -\dot R(0) R(0)/W .\label{tauchap}
\end{equation}
After  $\tau_0$, we choose the plus sign in \eqref{eqrho}, and we
find that
\begin{equation}
\label{solrho+}
\rho(\tau)W - E = (\tau- \tau_0)^2 W^2.
\end{equation}
Substituting the value of $\tau_0$ into \eqref{solrho+}, we find that
the general expression of the solution of \eqref{eqrho}  is given by
\begin{equation}
\label{solrho}
\rho(\tau)=\frac{E}{W} + W\bigl(\tau-\tau_0)^2.
\end{equation}
The angle $\theta$ is defined by \eqref{polar}; we let
\begin{equation}
\Theta(\tau)=\theta(t).\label{eqTheta}
\end{equation}
We are only interested for the present moment in the  principal part
of $\Theta$; it is a function $\Theta_1$ which satisfies the ordinary
differential equation:
\begin{equation}
\dot \Theta_1(\tau)=\frac{\sqrt E}{\rho},
\label{eqTheta0}\end{equation}
with the initial condition
\begin{equation}
\Theta_1(0)=0.\label{iniTheta0}
\end{equation}
We substitute the value of $\rho$ given by \eqref{solrho} into \eqref{eqTheta0}
and we find that
\begin{equation}
\dot \Theta_1(\tau)=\frac{W\sqrt E}{E +\bigl(W\tau +\dot R(0)
R(0)\bigr)^2}=\frac{\sqrt E}{R_1(\tau)^2}\label{eq:8}
\end{equation}
which we integrate immediately into
\begin{equation*}
\Theta_1(\tau)=\arctan\frac{W\tau + \dot R(0)
R(0)}{\sqrt E} - \arctan\frac{ \dot R(0)
R(0)}{\sqrt E}.
\end{equation*}
A more convenient way to write $\Theta_1$ is the following:
\begin{equation}
\Theta_1(\tau)=\arctan\frac{W(\tau-\tau_0)}{\sqrt E}
+\arctan\frac{W\tau_0}{\sqrt E}. \label{formTheta}
\end{equation}

Thanks to \eqref{initR} and \eqref{initdR}, we can see that 
\begin{equation}
W \sim \dot s(0)^2\eta^{-2}.\label{valE1.1}
\end{equation}
But \eqref{valE1.1} and \eqref{tauchap} imply that
\begin{equation}
\tau_0=O(\eta^4).\label{eq:20}
\end{equation}
This shows that we shall have to consider different cases: $\bar
\theta<\pi/2$ and $\bar \theta\ge\pi/2$.
If $\bar\theta<\pi/2$, we may suspect that $u_k$ will exit region
$\calr_2$ at time approximately $t_0+\bigl(\sqrt E \tan \bar
\theta\bigr)/\bigl(W\sqrt k\bigr)$; while if $\bar \theta\ge \pi/2$,
it is obvious that the ansatz is not sufficient: what will happen is
that $E/R^3$ is no more large with respect to $R$ and $\dot R$.

\subsection{Kernel for the linearized equation}\label{sec:kern-line-equat}

In order to go further, we have to validate asymptotics: consider
therefore the linear differential equation obtained from the
linearization of \eqref{ansatz}:
\begin{equation}
\ddot z + \frac{3Ez}{R_1^4} =0.\label{ansatzlin}
\end{equation}
We can obtain two linearly independent solutions by the following
argument: if we differentiate the ordinary differential equation
\eqref{ansatz} with respect to time, we find that $z_1=\dot R_1$ is a
solution of \eqref{ansatzlin}; we seek another solution of
\eqref{ansatzlin} under the form 
\begin{equation}\label{S2}
z_2= z_1 T.
\end{equation}
indeed, formula \eqref{solrho} gives the following form of $R_1$:
\begin{equation}
R_1(\tau)=\sqrt{\frac{E}{W}+ W(\tau -\tau_0)^2}.\label{solR}
\end{equation}
We differentiate this relation with respect to $\tau$, and we find
that
\begin{equation}
\dot R_1(\tau)=\frac{W(\tau -\tau_0)}{R_1}.\label{solR'}
\end{equation}
The equation satisfied by $T$ is
\begin{equation*}
2 \dot z_1 \dot T + z_1 \ddot T=0.
\end{equation*}
We multiply this equation by $z_1$, we integrate, and we find that, up
to an irrelevant multiplicative constant, $T$ satisfies the ordinary
differential equation:
\begin{equation*}
\dot T = \frac{1}{\dot R_1^2}. 
\end{equation*}
This can be integrated exactly and we obtain
\begin{equation*}
T=-\frac{E}{W^3(\tau -\tau_0)} + \frac{\tau}{W}.
\end{equation*}
According to definition \eqref{S2}, we have to multiply the above
expression by $\dot R_1$, for which we take expression \eqref{solR'}.
We obtain eventually
\begin{equation*}
z_2=-\frac{E}{W^2 R_1} + \frac{\tau(\tau - \tau_0)}{R_1}.
\end{equation*}
The Wronskian of $z_1$ and $z_2$ is readily calculated and is equal to
\begin{equation*}
z_1 \dot z_2 - z_2 \dot z_1 = z_1^2 \dot T =1;\label{wro}
\end{equation*}
therefore $z_1$ and $z_2$ are independent.
From here, we seek a kernel $K(\tau,\sigma)$ for $\tau\ge
\sigma$ which
satisfies the following conditions:
\begin{equation*}
\left\{
\begin{split}
&\frac{\partial^2K(\tau,\sigma)}{\partial \tau^2}
+\frac{3EK(\tau,\sigma)}{R_1^4(\tau)}  =0,\\
&K(\sigma,\sigma)=0,\\
&\frac{\partial K(\tau,\sigma)}{\partial \tau}\bigm|_{\tau=\sigma}=1,
\end{split}
\right.
\end{equation*}
under the form
\begin{equation*}
K(\tau,\sigma)=a_1(\sigma)S_1(\tau) + a_2(\sigma)S_2(\tau).
\end{equation*}
Thanks to \eqref{wro} and the definition of $z_1$,
we can see now that
\begin{equation*}
K(\tau,\sigma)=-z_2(\sigma) z_1(\tau) + z_1(\sigma) z_2(\tau),
\end{equation*}
which can be rewritten as
\begin{equation}
K(\tau,\sigma)=\frac{\tau-\sigma}{R_1(\tau)R_1(\sigma)}\left[\frac{E}{W}
+ W(\sigma -\tau_0)(\tau-\tau_0)\right], \quad \tau\ge
\sigma\ge 0.
\label{noyau}
\end{equation}
We extend $K(\tau,\sigma)$ by $0$ for $\sigma >\tau$.
It is convenient to define
\begin{equation}
J(\tau, \sigma)=\frac{E}{W} + W(\sigma-\tau_0)(\tau-\tau_0);\label{eq:1}
\end{equation}
with this notation, \eqref{noyau} becomes
\begin{equation}
K(\tau,\sigma)=\frac{(\tau-\sigma)J(\tau,\sigma)}{R_1(\tau)R_1(\sigma)},\label{eq:2}
\end{equation}
and we may also remark that
\begin{equation}
R_1(\tau)=\sqrt{J(\tau,\tau)}.\label{eq:3}
\end{equation}

\begin{remark}\label{Kge0}
The function $K(\tau, \sigma)\ge 0$ if $\tau\ll 1$:
indeed the only possibility for $K$ to be strictly negative is when the
product $(\sigma -\tau_0)(\tau-\tau_0)$ is strictly negative,
and $\sigma <\tau$; therefore, $\sigma$ is smaller than $\tau_0$ and $\tau$ is larger than $\tau_0$. Relations \eqref{tauchap}
and \eqref{valE1.1}
imply that $W^2(\tau_0-\sigma)(\tau -\tau_0)$ is estimated by
$C\tau$. Therefore, if $\tau\ll 1$, $K(\tau,\sigma)$ is nonnegative.
\end{remark}

\subsection{Validation of the earliest
asymptotic}\label{sec:valid-earl-asympt} 

With the help of the kernel $K$, we consider now the following problem:
to find an interval $[0, \tau_1]$ and a mapping $S_1$ from this
interval to $\Er$ such that $R_1+S_1$ solves \eqref{eqrnsc2}, with the
initial conditions \eqref{initR} and \eqref{initdR}. 
As $R_1(0)=R(0)$ and $\dot R_1(0)=\dot R(0)$,  $S_1(0)$ and $\dot
S_1(0)$ have to vanish; therefore \eqref{eqrnsc2}
can be rewritten as the following integral equation
\begin{equation}
S_1(\tau)=\call_1(S_1)(\tau) +G_1(\tau),\label{eqinteg}
\end{equation}
where $G_1$ is the function
\begin{equation}
G_1(\tau)=-\int_0^\tau K(\tau,\sigma)(2\alpha \dot R_1 +
R_1)(\sigma)\,d\sigma,
\label{call0}
\end{equation}
and $\call_1$ is an integral operator defined by
\begin{equation}
\begin{split}
\call_1(S_1)(\tau)&=\int_0^\tau \Bigl\{K(\tau,\sigma)\Bigl[\frac{3ES_1}{R_1^4}
-\frac{E}{R_1^3} +\frac{E(1-\ve)^2}{(R_1+S_1)^3}
- S_1
\Bigr](\sigma)\\&\qquad+ 2\alpha\frac{\partial
K(\tau,\sigma)}{\partial\sigma}S_1(\sigma) \Bigr\}\, d\sigma.
\end{split}\label{eq:15}
\end{equation}

Define
\begin{equation}
\tau_1 =\eta^{\gamma_1}, \quad 1 <{\gamma_1} < 2. \label{tauinf}
\end{equation}

Our purpose now is to prove that 
\eqref{eqinteg} has a unique solution on
$[0,\tau_1]$ thanks to the strict contraction principle.

We equip the space of continuous functions on $[0,\tau_1]$ with
the norm
\begin{equation}
\|S_1\|=\max\{|S_1(\sigma)|/R_1(\sigma): 0\le \sigma \le \tau_1\}.\label{norm}
\end{equation}
The choice of the weight $R_1$ in the norm is natural since we expect
$R_1$ to be the principal part of the solution; thus we expect that
the relative error $(R-R_1)/R_1$ will be small: our norm measures
precisely this relative error.

In order to apply the strict contraction principle, we estimate
certain functions through a sequence of technical calculations.

\begin{lemma}\label{borne} For all large enough $k$, the expression 
\begin{equation}
\label{cle}
I(\tau)=\frac{1}{R_1(\tau)}\int_0^ \tau \frac{K(\tau,\sigma)}{R_1^3(\sigma)} \, d\sigma.
\end{equation}
is bounded  on $[0,1]$. The bound will be called henceforth $\delta$.
\end{lemma}

\begin{remark}\label{rem:1} The expression $I(\tau)$ controls the
nonlinear term in the integral equation \eqref{eqinteg}, whose detail
is given in \eqref{eq:15}.
\end{remark}

\begin{proof}
We use the explicit expression of $K$ to perform an estimate of $I$:
\begin{equation*}
I(\tau)=\frac{1}{R_1^2(\tau)}\int_0^\tau
\frac{\tau-\sigma}{R_1^4(\sigma)} \left[
\frac{E}{W} + W(\sigma - \tau_0)(\tau -\tau_0)
\right]\,d\sigma.
\end{equation*}

We introduce the notation
\begin{equation}
\kappa=\sqrt E/W,\label{eq:4}
\end{equation}
and the change of variable
\begin{equation}
\tau=\tau_0 +\kappa y, \quad \sigma =\tau_0+\kappa x.\label{eq:5}
\end{equation}
The integral $I(\tau)$ is now given by
\begin{equation*}
\begin{split}
I(\tau)&=\frac{1}{J(\tau,\tau)} \int_0^\tau
\frac{(\tau-\sigma)J(\tau,\sigma)}{J(\sigma,\sigma)^2}\, d\sigma\\
&=\frac{1}{E} \frac{1}{1+y^2}\int_{-\tau_0/\kappa}^y
\frac{(y-x)(1+xy)}{(1+x^2)^2}\, dx.
\end{split}
\end{equation*}
We use the obvious inequalities
\begin{align}
&\abs{1+xy}\le\sqrt{1+x^2}\sqrt{1+y^2},\label{eq:6}\\
&\abs{y-x}\le \abs{x}+\abs{y},\label{eq:7}
\end{align}
to infer that
\begin{equation*}
\vabs{I(\tau)} \le\frac{1}{E} \frac{\abs{y}}{\sqrt{1+y^2}}
\int_{-\tau_0/\kappa}^y \frac{dx}{(1+x^2)^{3/2}}
+\frac{1}{E}\frac{1}{\sqrt{1+y^2}} \int_{-\tau_0/\kappa}^y
\frac{\abs{x}\, dx}{(1+x^2)^{3/2}}.
\end{equation*}
It is now clear that $\vabs{I(\tau)}$ is bounded independently of $y$,
i.e. of $\tau$ by a certain number $\delta$.
\end{proof}

We shall prove now that for  large enough
$k$, we can apply the strict contraction principle to 
 the mapping $S_1\mapsto \call_1(S_1) +G_1$ on the interval
$[0,\tau_1]$, defined by \eqref{tauinf}; for this purpose, we 
use the norm defined at \eqref{norm}, and we  show
the following result

\begin{theorem}\label{thr:1} For all $\gamma_1 \in (1,2)$,  and for
all small enough 
$p\in (0,1)$ 
there exists $k_0>0$
such
that for all $k\ge k_0$, the mapping $S_1\mapsto \call_1(S_1) +G_1$ leaves invariant the ball of
center $0$ and radius $p$ and is a strict contraction in that ball.
In particular, if $R$ is the solution of \eqref{eqrnsc2},
\eqref{initR} and \eqref{initdR}, we have the estimate
\begin{equation}
\max_{0\le \tau \le \eta^{\gamma_1}}
\vabs{\frac{R(\tau)-R_1(\tau)}{R_1(\tau)}}\le p.\label{eq:9}
\end{equation}
\end{theorem}

The proof of this result depends on several estimates given in
successive lemmas.

\begin{lemma} For all $\gamma_1 \in (1,2)$, there exists a constant $C$
such that
\begin{equation}
\|G_1\|\le C\bigl(\eta^{4\gamma_1 - 4} +\eta^{\gamma_1}\bigr).
\label{boule1}\end{equation}
\end{lemma}

\begin{proof}
By an integration by parts, 
\begin{equation*}
G_1 = 2\alpha K(\tau, 0) R_1(0) +\int_0^\tau \left(2\alpha
\frac{\partial K}{\partial \sigma}(\tau, \sigma) - K(\tau,
\sigma)\right) R_1(\sigma)\,d\sigma.
\end{equation*}
We estimate $|G_1(\tau)|/R_1(\tau)$: we first observe that
\begin{equation*}
\frac{K(\tau, 0)R_1(0)}{R_1(\tau)}
=\frac{\tau}{R_1^2(\tau)}\left[\frac{E}{W} - W\tau_0 (\tau
-\tau_0)
\right];
\end{equation*}
We estimate $R_1^2(\tau)$ from below by $E/W$; we also observe that
$W^2\tau_0$ is bounded, thanks to estimates \eqref{tauchap} and
\eqref{valE1.1}; therefore we have the following estimate, where we
have used again remark \ref{Kge0}:
\begin{equation}
\frac{K(\tau, 0)R_1(0)}{R_1(\tau)} \le \tau (1 + W^2\tau_0 \tau/E)
= 0(\eta^{\gamma_1}).\label{boule1.1}
\end{equation}

Next step is to calculate $\partial K(\tau, \sigma)/\partial \sigma$:
we use formulas \eqref{noyau}, \eqref{eq:1} and \eqref{eq:2} and we
find that 
\begin{equation*}
\frac{\partial K}{\partial \sigma}(\tau, \sigma)=L_1(\tau,\sigma)+
L_2(\tau,\sigma)+ L_3(\tau,\sigma),
\end{equation*}
where the $L_j$'s are respectively given by
\begin{align*}
L_1(\tau,\sigma) & =-\frac{J(\tau,\sigma)}{R_1(\tau)R_1(\sigma)},\\
L_2(\tau,\sigma)& =-\frac{(\tau - \sigma)(\sigma -\tau_0)
J(\tau,\sigma)W}{R_1(\tau)R_1^3(\sigma)},\\
L_3(\tau,\sigma)&=\frac{(\tau-\sigma)(\tau-\tau_0)W}{R_1(\tau)R_1(\sigma)}.
\end{align*}
Our aim is now to estimate the expressions
\begin{equation*}
I_j(\tau)=\frac{1}{R_1(\tau)} \int_0^\tau \vabs{L_j(\tau,\sigma)}
R_1(\sigma)\, d\sigma.
\end{equation*}
The first expression $I_1(\tau)$ is rewritten with the help of the
change of variable \eqref{eq:5} and becomes
\begin{equation*}
I_{1}(\tau)=\frac{\kappa}{1+y^2} \int_{-\tau_0/\kappa}^y \vabs{1+xy}\,
dx,
\end{equation*}
which we estimate as follows:
\begin{equation*}
I_1(\tau)\le \frac{\kappa}{\sqrt{1+y^2}} \int_{-\tau_0/\kappa}^y
\sqrt{1+x^2}\, dx.
\end{equation*}
Since $\tau_0/\kappa\le 1$ for $k$ large enough, we can see that
\begin{equation*}
\begin{split}
\frac{1}{\sqrt{1+y^2}}\int_{-\tau_0/\kappa}^y \sqrt{1+x^2}\, dx
&\le \sqrt 2 \int_{-\tau_0/\kappa}^{\min(y,0)} \,
dx+\int_{\min(y,0)}^y \frac{\sqrt{1+x^2}}{\sqrt{1+y^2}}\, dx\\ 
&\le \sqrt 2\left(\frac{\tau_0}{\kappa} +y\right).
\end{split}
\end{equation*}
Therefore
\begin{equation*}
I_1(\tau)\le \sqrt 2\, \tau.
\end{equation*}

With the change of variable \eqref{eq:5}, we may write the expression
$I_2(\tau)$ as
\begin{equation*}
I_2(\tau)=\frac{\kappa}{(1+y^2)}\int_{-\tau_0/\kappa}^y
\frac{(y-x)\abs{x}\abs{1+xy}\, dx}{1+x^2}.
\end{equation*}
We use \eqref{eq:6} again together with
\begin{equation*}
\frac{\abs{x}}{\sqrt{1+x^2}} \le 1,
\end{equation*}
and we find
\begin{equation*}
I_2(\tau)\le \frac{\kappa}{\sqrt{1+y^2}}\int_{-\tau_0/\kappa}^y
(y-x)\, dx \le \frac{W\tau^2}{2\sqrt E\,
\sqrt{1+(\tau-\tau_0)^2/\kappa^2}}. 
\end{equation*}
But for $k$ large enough and $\tau\ge 0$,
\begin{equation*}
\sqrt{1+(\tau-\tau_0)^2/\kappa^2} \ge\frac{\tau}{2\kappa};
\end{equation*}
indeed, if $\tau\ge 2\tau_0$, $\tau-\tau_0\ge \tau/2$, and the
inequality is clear; on the other hand, if $\tau\le 2\tau_0$, for $k$
large enough $\tau_0\le \kappa$, and the inequality also follows.
Therefore, there exists a number $C$ such that for all large enough
$k$ and all $\tau$ in $[0,\eta^{\gamma_1}]$ the following inequality
holds:
\begin{equation*}
I_2(\tau)\le C\eta^{\gamma_1}.
\end{equation*}

The third expression is handled as follows:
\begin{equation*}
I_3(\tau)=\frac{1}{R_1(\tau)} \int_0^\tau
\frac{(\tau-\sigma)\abs{\tau-\tau_0}W \, d\sigma}{R_1(\tau)}
=\frac{\abs{\tau-\tau_0}\tau^2 W}{2J(\tau,\tau)}. 
\end{equation*}
If $0\le \tau\le \tau_0 +\kappa$, we use the inequality
$J(\tau,\tau)\ge E/W$ and we obtain
\begin{equation*}
I_3(\tau)\le \frac{\kappa^2 \tau^2 W^2}{2E} =O(\eta^2),
\end{equation*}
since $\abs{\tau-\tau_0}\le \max(\kappa,\tau_0)=\kappa$ for $k$ large
enough. On the other hand, for $\tau \ge \tau_0+\kappa$ and for $k$
large enough
\begin{equation*}
\frac{\tau}{\tau-\tau_0} \le 2,
\end{equation*}
and therefore, using the inequality $J(\tau,\tau)\ge W\abs{\tau
-\tau_0}^2$, we obtain
\begin{equation*}
I_3(\tau) \le\frac{\abs{\tau-\tau_0}\tau^2 W}{2\abs{\tau-\tau_0}^2
W}\le \tau.
\end{equation*}
Thus, we have shown that
\begin{equation}
\frac{1}{R_1(\tau)}\int_0^\tau \vabs{\frac{\partial K}{\partial
\sigma}(\tau,\sigma)}R_1(\sigma) \, d\sigma =0\bigl(\eta^{\gamma_1}\bigr).\label{eq:10}
\end{equation}

There remains to estimate
\begin{equation}
\frac{1}{R_1(\tau)}\int_0^\tau K(\tau,\sigma)R_1(\sigma)\, d\sigma.\label{boule''}
\end{equation}
We rewrite \eqref{boule''} as
\begin{equation*}
\frac{1}{R_1^2(\tau)}\int_0^\tau R_1^4(\sigma) \frac{\tau
-\sigma}{R_1^4(\sigma)}
\left[\frac{E}{W} +
 W(\sigma-\tau_0)(\tau -\tau_0)\right]\, d\sigma
\end{equation*}
and we find that thanks to remark \ref{Kge0}
\begin{equation}
\begin{split}
&\frac{1}{R_1(\tau)} \int_0^\tau \abs{K(\tau,\sigma)} R_1(\sigma)\,
d\sigma\\
&\quad=\frac{1}{R_1(\tau)} \int_0^\tau K(\tau,\sigma)R_1(\sigma) \le
\delta\max(R_1^4(0), R_1^4(\tau)) =O\bigl(\eta^{4\gamma_1 -4}\bigr).
\end{split}
\label{eq:11}
\end{equation}

Summarizing \eqref{eq:11} with \eqref{boule1.1} and \eqref{eq:10},
we find estimate \eqref{boule1}.
\end{proof}

Next lemma enables us to estimate $\|\call_1(S_1)\|$ when $\|S_1\|\le p<1$.

\begin{lemma} Assume $\norm{S_1}\le p<1$. For all $\gamma_1 \in
(1,2)$, there exists $k_0$ and $C$ such that 
for all $k\ge k_0$ the following estimate holds:
\begin{equation*}
\bigl\|\call_1(S_1)\bigr\| \le \frac{6E\delta p^2}{(1-p)^5} +\frac {2\ve
E\delta}{(1-p)^3} + pC\bigl(\eta^{4\gamma_1 -4}+\eta^{\gamma_1}\bigr)
\end{equation*}
\end{lemma}

\begin{proof}
The easiest part is the estimate on 
\begin{equation*}
\int_0^\tau \left(-K(\tau, \sigma) S_1(\sigma) + 2\alpha
\frac{\partial K}{\partial \sigma}(\tau,\sigma)S_1(\sigma)\right)\, d\sigma.
\end{equation*}
We can see that the absolute value of this expression is estimated by
\begin{equation*}
\|S_1\|\int_0^\tau \left(2\alpha
\left|\frac{\partial K}{\partial \sigma}(\tau,\sigma)\right| +
K(\tau,\sigma)\right)R_1(\sigma)\, d\sigma .
\end{equation*}
We recognize expressions which have already been estimated in
\eqref{eq:10} and \eqref{eq:11}. Therefore, it
is immediate that
\begin{equation}
\label{boule2.1}
\begin{split}
&\left|\frac{1}{R_1(\tau)}\int_0^\tau 
\left(-K(\tau, \sigma) S_1(\sigma) + 2\alpha
\frac{\partial K}{\partial \sigma}(\tau,\sigma)S_1(\sigma)\right)
\, d\sigma
\right|\\
&\quad\le C\|S_1\|\bigl(\eta^{4\gamma_1 -4}+\eta^{\gamma_1}\bigr).
\end{split}
\end{equation}

Next comes the slightly more complicated expression
\begin{equation}
\frac{1}{R_1(\tau)}\int_0^\tau K(\tau, \sigma) \frac{(2\ve
-\ve^2)E}{(R_1+S_1)^3(\sigma)}\, d\sigma. \label{compli1}
\end{equation}

For all $k>0$, $2\ve-\ve^2$ is at most equal to $2\ve$; therefore, if
$\|S_1\|\le p <1$, then the absolute value of \eqref{compli1}
is estimated by
\begin{equation*}
\frac{2\ve  E}{R_1(\tau)}\int_0^\tau \frac{K(\tau,
\sigma)}{R_1^3(\sigma)(1-p)^3}\, d\sigma.
\end{equation*}

Therefore, we can see that
\begin{equation}
\left|\frac{1}{R_1(\tau)}\int_0^\tau K(\tau, \sigma) \frac{(2\ve
-\ve^2)E}{(R_1+S_1)^3(\sigma)}\, d\sigma
\right|\le \frac{2\ve  E \delta}{(1-p)^3}.
\label{boule2.2}
\end{equation}

The last and most complicated term contains the expression
\begin{equation*}
\left(\frac{E}{(R_1+S_1)^3}-\frac{E}{R_1^3} +
\frac{3ES_1}{R_1^4}\right)(\sigma),
\end{equation*}
which can be rewritten thanks to Taylor's formula with integral
remainder as
\begin{equation*}
12E\int_0^1
 \frac{S_1^2(\sigma)}{\bigl(R_1(\sigma)+sS_1(\sigma)\bigr)^5}(1-s)\, ds.
\end{equation*}
The absolute value of this expression is estimated by
\begin{equation*}
\frac{6Ep^2}{(1-p)^5 R_1(\sigma)^3}.
\end{equation*}
Therefore, we find that
\begin{equation}
\left| \frac{1}{R_1(\tau)}\int_0^\tau K(\tau, \sigma)\left(\frac{E}{(R_1+S_1)^3}-\frac{E}{R_1^3} +
\frac{3ES_1}{R_1^4}\right)(\sigma)\, d\sigma
\right|
\le \frac{6E\delta p^2}{\bigl(1-p\bigr)^5}.\label{boule2.3}
\end{equation}

\end{proof}

\begin{proof}[End of the proof of Theorem \ref{thr:1}]

The ball of radius $p$ about $0$ will be invariant by the mapping
$S_1\mapsto \call_1(S_1)+G_1$ provided that
\begin{equation}
\begin{split}
&\frac{6E\delta p^2}{\bigl(1-p\bigr)^5} +  \frac{2\ve  E
\delta}{(1-p)^3}
+ pC \bigl(\eta^{4\gamma_1 -4}+\eta^{\gamma_1}\bigr) \\
&\quad +
C\bigl(\eta^{4\gamma_1 -4}+\eta^{\gamma_1}\bigr)
\le p.
\end{split}\label{invarboule}
\end{equation}
Choose $p$ small enough for the following inequality to hold
\begin{equation*}
\frac{6E\delta p}{\bigl(1-p\bigr)^5} \le \frac{1}{2};
\end{equation*}
Choose then $k$ so large that \eqref{invarboule} holds.

Let us prove now that for an adequate choice of $p$ and $k_0$ and for
all $k\ge k_0$, the mapping  $\call_1$ is a strict contraction: the
easiest part of the estimate pertains to
\begin{equation*}
\frac{1}{R_1(\tau)} 
\int_0^\tau 
\left(
2\alpha (S_1-\hat S_1)(\sigma)
\frac{\partial K}{\partial\sigma }
(\tau,\sigma) - (S_1-\hat S_1)(\sigma)K(\tau,\sigma)
\right)\, d\sigma
\end{equation*}
and it is clear from the proof of estimate \eqref{boule2.1} that
\begin{equation}
\begin{split}
&\frac{1}{R_1(\tau)} \left|\int_0^\tau 
\left(
2\alpha (S_1-\hat S_1)(\sigma)
\frac{\partial K}{\partial\sigma }
(\tau,\sigma) - (S_1-\hat S_1)(\sigma)K(\tau,\sigma)
\right)
\, d\sigma \right|\\
&\quad \le C\|S_1 -
\hat S_1\|\bigl(\eta^{4\gamma_1 -4}+ \eta^{\gamma_1}\bigr).
\end{split}
\label{boule3.1}
\end{equation}

The second easiest term involves the difference
\begin{equation*}
\frac{2\ve -\ve^2}{(R_1+S_1)^3} -\frac{2\ve -\ve^2}{(R_1+\hat S_1)^3}
\end{equation*}
which we estimate thanks to Taylor's formula:
\begin{equation*}
6\ve \int_0^1 \frac{|(\hat S_1 - S_1)(\sigma)|}{(R_1(\sigma) +
S_1(\sigma) +s(\hat S_1(\sigma) - 
S_1(\sigma)))^4}\, ds \le \frac{6\ve \|\hat S_1 - S_1\|}{(1-p)^4
R_1(\sigma)^3}. 
\end{equation*}
Therefore, the corresponding term in $\call_1(S_1)-\call_1(\hat S_1)$
contributes an estimate
\begin{equation}
\frac{6\ve E\delta\|\hat S_1 - S_1\|}{(1-p)^4 }.\label{boule3.2}
\end{equation}
The last and most complicated term involves the expression
\begin{equation*}
\frac{3\hat S_1}{R_1^4} +\frac{1}{(R_1+\hat S_1)^3} -\frac{3S_1}{R_1^4} -
\frac{1}{(R_1+S_1)^3}.
\end{equation*}
We use a Taylor expansion twice to rewrite this expression as
\begin{equation*}
12 (\hat S_1 - S_1)\int_0^1 \int_0^1 \frac{\bigl(S_1 + s(\hat S_1 - S_1)\bigr)\, ds'\, ds}
{\bigl(R_1 +s'(S_1+s(\hat S_1-S_1))\bigr)^5}.
\end{equation*}
Therefore, the corresponding term is estimated by 
\begin{equation*}
\frac{12p\|\hat S_1 - S_1\|}{(1-p)^5 R_1(\sigma)^3}.
\end{equation*}
Thus, the norm of the corresponding contribution in $\call_1(\hat S_1)
-\call_1(S_1)$ is estimated by
\begin{equation}
\frac{12p E\delta\|\hat S_1 - S_1\|}{(1-p)^5 }.\label{boule3.3}
\end{equation}

Therefore, if we summarize the estimates \eqref{boule3.1},
\eqref{boule3.2} and \eqref{boule3.3}, we find that on a ball of
radius $p<1$ about $0$, the Lipschitz constant of $\call_1$ is estimated
by
\begin{equation}
C\bigl(\eta^{4\gamma_1 -4}+ \eta^{\gamma_1}\bigr) + \frac{6\ve
E\delta}{(1-p)^4 }  
+\frac{12p E\delta}{(1-p)^5 }.\label{boule3}
\end{equation}
If we choose $p$ small enough for $12p E\delta/(1-p)^5$ to be less
than or equal to $1/2$, it is clear that we can choose $k_0$ large
enough for the sum of the remaining terms in \eqref{boule3} to be less
than or equal to $1/4$. 

Together with the conditions found above for the invariance of the
ball of radius $p$ about $0$, we have shown the first part of theorem
\ref{thr:1}. We also infer from this proof that
\begin{equation}
\vnorm{S_1} \le C\bigl(\eta^{4\gamma_1-4} +\eta^{\gamma_1}+\ve\bigr).\label{eq:23}
\end{equation}
Its last assertion is an immediate consequence of the equivalence of
\eqref{eqinteg} with \eqref{eqrnsc2}, the definition of the norm, and
the fact that the initial data coincide.
\end{proof}

\begin{remark}\label{ppetit}
Let us observe that in the end of the proof of Theorem
\ref{thr:1}, we can take $p$ arbitrarily small provided that
$k_0$ is chosen large enough.
\end{remark}

We conclude this section by the

\begin{proposition}\label{thr:2} Let $\gamma_1$ belong to $(1,2)$, and let
$\tau_1=\eta^{\gamma_1}$. Then the following equivalences hold:
\begin{align}
R(\tau)&\sim R_1(\tau) \text{ over $[0,\tau_1]$,} \label{eq:21}\\
\dot R(\tau) &\sim \dot R_1(\tau) \text{ over $[\eta^3,\tau_1]$.}\label{eq:22}
\end{align}
\end{proposition}

\begin{proof}
The first statement is an almost immediate consequence of
\eqref{eq:23}: we have
\begin{equation*}
\frac{R(\tau)}{R_1(\tau)} = 1+\frac{S_1(\tau)}{R_1(\tau)},
\end{equation*}
so that on $[0,\tau_1]$
\begin{equation*}
\vabs{\frac{R(\tau)}{R_1(\tau)}-1} \le \vnorm{S_1} \le
C\bigl(\eta^{4\gamma_1-4}+\eta^{\gamma_1}+\ve\bigr), 
\end{equation*}
and \eqref{eq:21} follows.

In order to compare $\dot R$ and $\dot R_1$, we write the differential equations
that they satisfy:
\begin{align*}
\ddot R + 2\alpha \dot R &= - R +\frac{E(1-\ve)^2}{R^3},\\
\ddot R_1 + 2\alpha \dot R_1&=\frac{E}{R_1^3} + 2\alpha \dot R_1.
\end{align*}
Therefore, if we subtract the second of these equations from the
first, we deduce that
\begin{equation*}
\ddot R-\ddot R_1 + 2\alpha(\dot R-\dot R_1) = - R +\frac{E(1-\ve)^2}{R^3}
-\frac{E}{R^3_1} - 2\alpha \dot R_1.
\end{equation*}
The initial data vanish.

We integrate because we want to estimate $\dot R - \dot R_1$:
\begin{equation}
\begin{split}
&\dot R(\tau) -\dot R_1(\tau)\\
&\quad=
\int_0^{\tau}
\exp\bigl(-2\alpha(\tau-\sigma)\bigr)
\Bigl\{
 - R
\quad+\frac{E(1-\ve)^2}{R^3}
-\frac{E}{R^3_1}
- 2\alpha \dot R_1
\Bigr\}(\sigma)\, d\sigma.
\end{split}\label{deriv}
\end{equation}

The next step is to estimate the integral on the right hand side of
\eqref{deriv}. We decompose this integral into three terms:
\begin{align*}
I_1&=-2\alpha\int_0^{\tau}\exp\bigl(-2\alpha(\tau-\sigma)\bigr) \dot
R_1(\sigma)\, d\sigma,\\
I_2&=-\int_0^{\tau}\exp\bigl(-2\alpha(\tau-\sigma)\bigr)
R(\sigma)\, d\sigma,\\
I_3&=E\int_0^{\tau}\exp\bigl(-2\alpha(\tau-\sigma)\bigr)
\left(\frac{(1-\ve)^2}{R^3} - \frac{1}{R_1^3}\right)(\sigma)\, d\sigma.
\end{align*}

The first two integrals are very easy to estimate: for $I_1$ an
integration by parts gives
\begin{equation*}
I_1=-2\alpha R_1(\tau) + 2\alpha e^{-2\alpha \tau} R_1(0)
-4\alpha^2\int_0^{\tau} \exp\bigl(-2\alpha(\tau-\sigma)\bigr)
R_1(\sigma)\, d\sigma.
\end{equation*}
Thanks to \eqref{solR}, on $[0,\tau_1]$,
$R_1(\tau)=O\bigl(\eta^{\gamma_1-1}\bigr)$. Therefore
\begin{equation}
\abs{I_1}= O\bigl(\eta^{\gamma_1-1}\bigr).\label{i1}
\end{equation}
For $I_2$, the situation is even simpler since it can be readily seen
that
\begin{equation}
\vabs{I_2}\le (1+p)\int_0^{\tau}
R_1(\sigma)\exp\bigl(-2\alpha(\tau-\sigma)\bigr)\, d\sigma
=O\bigl(\eta^{2\gamma_1-1}\bigr),\label{i2} 
\end{equation}
where $p=\norm{S_1}$.

There remains to estimate $I_3$; a straightforward calculation shows
that
\begin{equation*}
\vabs{\frac{(1-\ve)^2}{R^3} - \frac{1}{R_1^3}} \le \frac{2\ve 
+p(3+3p +p^2)}{(1-p)^3} \frac{1}{R_1^3};
\end{equation*}
here $p=\vnorm{S_1}$ is estimated at \eqref{eq:23}, so that
\begin{equation}
\vabs{\frac{(1-\ve)^2}{R^3} - \frac{1}{R_1^3}} =O\bigl(\ve
+\eta^{\gamma_1}+\eta^{4(\gamma_1-1)}\bigr). \label{i3}
\end{equation}
Therefore,
\begin{equation*}
I_3= O\bigl(\ve +\eta^{\gamma_1}+\eta^{4(\gamma_1-1)}\bigr) J,
\end{equation*}
where $J$ is defined as
\begin{equation*}
J=\int_0^{
\tau}\exp\bigl(-2\alpha(\tau-\sigma)\bigr)\frac{E}{R_1^3(\sigma)}
\,d\sigma.
\end{equation*}
But $J$ can be more conveniently rewritten as
\begin{equation*}
J=\int_0^{\tau}\exp\bigl(-2\alpha(\tau-\sigma)\bigr) \ddot R_1(\sigma)\, d\sigma,
\end{equation*}
which we integrate by parts. We find that
\begin{equation*}
\begin{split}
&J=\dot R_1(\tau) - e^{-2\alpha\tau} \dot R_1(0) +2\alpha\bigl(
R_1(\tau) - e^{-2\alpha\tau} R_1(0)\bigr)\\
&\qquad+4\alpha^2\int_0^{\tau} \exp\bigl(-2\alpha(\tau-\sigma)\bigr)
R_1(\sigma)\, d\sigma.
\end{split}
\end{equation*}
We can see now that
\begin{equation*}
J=\dot R_1(\tau)+O\bigl(\eta^{\gamma_1 -1}\bigr).
\end{equation*}
Therefore
\begin{equation*}
\frac{\dot R(\tau)}{\dot R_1(\tau)}
=\left(1+\frac{O\bigl(\eta^{\gamma_1-1}\bigr)}{\dot R_1(\tau)}\right)
\bigl(1 + O\bigl(\ve
+\eta^{\gamma_1}+\eta^{4(\gamma_1-1)}\bigr)\bigr). 
\end{equation*}

Since $\ddot R_1$ is nonnegative, for $\tau\in[\eta^3, \tau_1]$, $\dot
R_1(\tau)$ can be estimated from below by $\dot R_1\bigl(\eta^3\bigr)$
which is equal to
\begin{equation*}
\frac{W(\eta^3 -\tau_0)}{\sqrt{(E/W)+W(\eta^3 -\tau_0)^2}};
\end{equation*}
we infer from relations \eqref{valE1.1} and \eqref{eq:20} that
\begin{equation*}
\lim_{k\to \infty}\dot R_1\bigl(\eta^3\bigr)>0,
\end{equation*}
which concludes the proof.
\end{proof}

\section{The case $\bar \theta<\pi/2$.}\label{sec:case-bar-thetapi2}

We prove here the first theorem which justifies Moreau's rule for
$\bar \theta <\pi/2$.

\begin{theorem}
If $\bar \theta<\pi/2$, the representative point of the system enters
region $\calr_3$ at a time $\bar t=t_0 +\bar \tau/\sqrt k$, where 
\begin{equation}
\bar \tau \sim \tan\bar\theta\, \sqrt E/W; \label{eq:12}
\end{equation}
moreover, in the coordinates defined by the axes $y_1$ and $y_2$ (see
Fig. \ref{fig:region}), we have the following asymptotics  for all
$t\ge \bar t$:
\begin{align}
&0\le y_1(t)\le \frac{C\exp\bigl(\xi_1(t-\bar t)\sqrt k\bigr)}{\sqrt
k},\label{eq:13}\\ 
&y_2(t)\sim (t-\bar t)\dot s(0)\cos \bar \theta. \label{eq:14}
\end{align}
\end{theorem}

\begin{proof}
The differential equation satisfied by $\Theta$ defined by
\eqref{eqTheta} is deduced from \eqref{invar} and is given by
\begin{equation}
\dot \Theta =\frac{\sqrt E (1-\ve)}{R^2}, \quad \Theta(0)=0.\label{eq:33}
\end{equation}
Recall that the principal part $\Theta_1$ is defined by
\eqref{eqTheta0} and \eqref{iniTheta0}.
Let us estimate $\Upsilon_1 =\Theta- \Theta_1$: $\Upsilon_1$ satisfies the
differential equation
\begin{equation*}
\dot \Upsilon_1 = \frac{\sqrt E(1-\ve)}{R^2}-\frac{\sqrt E}{R_1^2}.
\end{equation*}
Therefore, if we let $p=\vnorm{S_1}$, and if we denote
\begin{equation}
\beta=\frac{\ve + p^2 + 2p}{(1-p)^2},\label{coin:eq:1}
\end{equation}
we find that
\begin{equation*}
\vabs{\dot\Upsilon_1}\le \beta \dot \Theta_1.
\end{equation*}
Hence, for all $\tau\in[0,\tau_1]$,
\begin{equation}
(1-\beta)\Theta_1(\tau) \le \Theta(\tau) \le (1+\beta) \Theta_1(\tau).\label{cadre}
\end{equation}

According to \eqref{eq:23} and the definition of $\ve$, there
exists $k_0$ such that for all $k\ge k_0$:
\begin{equation}
\beta<1, \quad \bar\theta <\frac{\pi(1-\beta)}{2}. \label{trucmuche}
\end{equation}

Let $\tau_+$ and $\tau_-$ be defined by the relations
\begin{equation*}
\Theta_1(\tau_+) = \frac{\bar\theta}{1-\beta}, \quad \Theta_1(\tau_-)
= \frac{\bar\theta}{1+\beta}.
\end{equation*}
Thanks to condition \eqref{trucmuche} and formula \eqref{formTheta}
$\tau_+$ and $\tau_-$ are well defined, and are given by
\begin{align*}
\tau_+ &= \tau_0 +\frac{\sqrt E}{W}
\frac{\tan\bigl(\bar\theta/(1-\beta)\bigr)-W\tau_0 /\sqrt
E}{1+(W\tau_0/\sqrt E)
\tan\bigl(\bar\theta/(1-\beta)\bigr) },\\
\tau_- &= \tau_0 +\frac{\sqrt E}{W}
\frac{\tan\bigl(\bar\theta/(1+\beta)\bigr)-W\tau_0 /\sqrt
E}{1+(W\tau_0/\sqrt E)
\tan\bigl(\bar\theta/(1+\beta)\bigr) }.
\end{align*}
Therefore, as $k$ tends to infinity, both $\tau_-$ and $\tau_+$ are
equivalent to $\bigl(\sqrt E / W)\tan \bar \theta$.

The function $\Theta$ is strictly increasing with respect to time;
thanks to inequality \eqref{cadre}, there is a unique $\bar \tau\in
[\tau_-,\tau_+]$
such that $\Theta(\bar \tau)=\bar \theta$. 

We know an equivalent of $\tau_-$ and $\tau_+$ as $k$ tends to
infinity:
\begin{equation*}
\bar \tau \sim \frac{\sqrt E}{W} \tan \bar \theta=O\bigl(\eta^2\bigr).
\end{equation*}

Together with \eqref{solR}, the above relation implies
\begin{equation}
R_1(\bar\tau) \sim \frac{\dot r(0)\eta}{2\cos\bar\theta
\sqrt{\alpha^2-1}}\label{valbar}
\end{equation}
and from \eqref{solR'} that 
\begin{equation}
\dot R_1(\bar \tau) \sim \frac{\dot s(0)\sin \bar \theta}{\eta}.\label{valdbar}
\end{equation}

Proposition \ref{thr:2} implies $R(\bar\tau)\sim R_1(\bar \tau)$ and
$\dot R(\bar \tau)\sim \dot R_1(\bar \tau)$.

We change coordinates now, taking the axis $y_2$ along the second side of the
convex cone $K$ and the axis $y_1$ perpendicular to $y_2$, and going
out of $K$. The new time variable is a translation of the natural
time, denoted $t'$, and we set its origin at
the time when the representative point enters region $\calr_3$.
We also let $\bar t= t_0+\bar \tau/\sqrt k$.

With these conventions, 
\begin{equation*}
y(0)= r(\bar t),\quad \dot y(0) = \dot r(\bar t) + i r(\bar t)
\dot\theta(\bar t).
\end{equation*}

We use now the equivalents obtained previously:
\begin{align}
&y_1(0)=O(\eta^2/\sqrt k),\quad y_2(0)=0,\label{coin:eq:13}\\
&\dot y_1(0)\sim \dot s(0)\sin\bar\theta,\quad \dot y_2(0)\sim
\dot s(0) \cos\bar\theta.\label{coin:eq:14}
\end{align}

The second component $y_2$ of $y$ satisfies the ordinary differential
equation
\begin{equation*}
\ddot y_2=0,
\end{equation*}
so that
\begin{equation}
y_2(t')\sim t' \dot s(0) \cos\bar\theta.\label{coin:eq:2}
\end{equation}
The first component $y_1$ of $y$ satisfies the following ordinary
differential equation
\begin{equation}
\ddot y_1 + 2 \alpha \sqrt k \dot y_1 + k y_1=0,\label{coin:eq:15}
\end{equation}
as long as $y_1\ge0$. The explicit solution of \eqref{coin:eq:15} with
initial data \eqref{coin:eq:13} and \eqref{coin:eq:14} is given by
\begin{equation*}
y_1(t')=\dot y_1(0)\frac{e^{\xi_1 t'\sqrt k}-e^{\xi_2t'\sqrt
k}}{2\sqrt k \sqrt{\alpha^2 -1}} +y_1(0)\frac{\xi_1 e^{\xi_2 t'\sqrt
k}- \xi_2 e^{\xi_1t'\sqrt
k}}{2 \sqrt{\alpha^2 -1}} .
\end{equation*}
Since $\dot y_1(0)$ is non negative, $y_1(t')$ stays non negative for
all $t'\ge 0$ and we have the following estimate on the first component of $y$:
\begin{equation}
0\le y_1(t') \le \frac{Ce^{-\abs{\xi_1}t'\sqrt k}}{\sqrt k}.
\label{coin:eq:3}
\end{equation}
\end{proof}

Thus, we obtain the conclusion of this section as the following
Theorem:

\begin{theorem}
Let $\Pi_1$ be the orthogonal projection on $\{x_1=0\}$, and let
$\Pi_2$ be the orthogonal projection on $\{x_1\cos \bar\theta + x_2
\sin \bar \theta =0\}$; see Fig. \ref{fig:region}. As $k$ tends to
infinity, $u_k$ converges 
uniformly on compact sets of $\Er^+$ to $u_\infty$ given by
\begin{equation*}
u_\infty(t)=\begin{cases}
u(0)+t\Pi_1 \dot u(0) &\text{if $0\le t \le t_0$,}\\
(t-t_0)\Pi_2 \Pi_1 \dot u(t_0)&\text{if $t_0\le t$.}
\end{cases}
\end{equation*}
\end{theorem}

\begin{proof}The initial part of the motion is described thanks to
\eqref{eqr} and \eqref{eqs}.
Estimate \eqref{coin:eq:3} proves that $y_1(t)$ tends to $0$ uniformly
on compact sets of $]t_0, \infty)$; relation
\eqref{coin:eq:2} enables us to conclude.
\end{proof}

\section{The second asymptotics}\label{sec:second}

In the case $\bar\theta\ge \pi/2$, we need a new asymptotic, and an
estimate which is based essentially on the use of Lyapunov
functionals, and which will be proved in Section
\ref{sec:final-asymptotics}. 

We restrict the choice of the exponent $\gamma_1$ in the definition of
$\tau_1$ by assuming that
\begin{equation}
\gamma_1\in \left(1,\frac{4}{3}\right).\label{gammaetoile}
\end{equation}
The reason for this choice is the following: if \eqref{gammaetoile}
holds, then the term $E(1-\ve)^2/R(\tau_1)^3$ is of order
$\eta^{3-3\gamma_1}$ which is small relatively to $\dot R(\tau_1)$,
according to the analysis of \ref{thr:2}: indeed, the following equivalents of
$R(\tau_1)$ and $\dot R(\tau_1)$ are a consequence of proposition
\ref{thr:2}:
\begin{equation}
\label{equiv*}
R(\tau_1)\sim \dot s(0) \eta^{\gamma_1-1}, \quad \dot
R(\tau_1)\sim  \dot s(0) \eta^{-1}.
\end{equation}

Let $\zeta$ be such that
\begin{equation}
0<\zeta<1/|\xi_1|.\label{Kun}
\end{equation}

We define the time $\tau_3$ by
\begin{equation}
\tau_3=\zeta \ln(1/\eta). \label{tau1}
\end{equation}
We use the notation $\tau_3$, because we will define below an
intermediate time $\tau_2$ between $\tau_1$ and $\tau_3$.

We claim that the solution of \eqref{eqrnsc2} on the interval
$[\tau_1,\tau_3]$ is very close to the solution of
\begin{equation*}
\ddot R_2 + 2\alpha \dot R_2 + R_2=0, \quad
R_2(\tau_1)= R(\tau_1),\quad
\dot R_2(\tau_1)=\dot R(\tau_1).
\end{equation*}

Let us define two kernels $K_2$ and $H_2$ on $\Er^+$ by
\begin{align*}
K_2(\tau)&=\frac{e^{\xi_1\tau}-e^{\xi_2\tau}}{2\sqrt \Delta},\\
H_2(\tau)&=\frac{-\xi_2 e^{\xi_1\tau}+\xi_1 e^{\xi_2\tau}}{2\sqrt \Delta}.
\end{align*}
We extend $K_2$ and $H_2$ to $\Er^-$ by $0$.
Therefore, $R_2$ is given explicitly for $\tau\ge \tau_1$ by
\begin{equation}
\label{R1}
R_2(\tau)=K_2(\tau-\tau_1)\dot R(\tau_1)+H_2(\tau-\tau_1) R(\tau_1).
\end{equation}

In order to substantiate our claim, we argue as for
theorem \ref{thr:1}: write $R=R_2 + S_2$; then $S_2$ is a solution of
the integral equation
\begin{equation}
\label{intR2}
S_2(\tau)=\int_{\tau_1}^\tau K_2(\tau-\sigma) \frac{E(1-\ve)^2}{(R_2+
S_2)^3(\sigma)}\, d\sigma.
\end{equation}

It is convenient to denote
\begin{equation*}
\call_2(S_2)=\int_{\tau_1}^\tau K_2(\tau-\sigma) \frac{E(1-\ve)^2}{(R_2+
S_2)^3(\sigma)}\, d\sigma,
\end{equation*}
whenever $R_2+S_2$ does not vanish over $[\tau_1,\tau_3]$.

Let us prove that $R_2$ never vanishes over $[\tau_1,+\infty)$:
thanks to the inequalities $0>\xi_1>\xi_2$, the functions
$K_2$ and $H_2$ are positive for $\tau>0$, $H_2(0)$ is equal to $1$; $R(\tau_1)$
and $\dot R(\tau_1)$ are strictly positive. Thus the positivity of $R_2$ is clear.

On the space $C^0([\tau_1,\tau_3])$, we introduce the norm
\begin{equation}
\norm{S_2}=\sup\{\abs{S_2(\tau)}/R_2(\tau): \tau\in
[\tau_1,\tau_3]\}.\label{norm1}
\end{equation}

We remark that $\call_2$ is well defined on the open ball of radius
$1$ about $0$ in the norm \eqref{norm1}.

We prove that $\call_2$ is a contraction on an appropriate ball, which
will lead us to validated asymptotics for $R$ on the interval
$[\tau_1,\tau_3]$. 

\begin{theorem}
\label{contractionbis}
For all $p\in (0,1)$, there exists $k_1>0$ such that for all $k>k_1$,
$\call_2$ is a contraction from the ball of radius $p$ (relatively to
$\norm{\>}$) about $0$ to itself.
\end{theorem}

\begin{proof}
We will show in Lemma \ref{thr:3} that the expression
\begin{equation*}
I(\tau)=\frac{1}{R_2(\tau)}\int_{\tau_1}^{\tau}
\frac{K_2(\tau-\sigma)}{R_2^3(\sigma)} \,d\sigma
\end{equation*}
tends to $0$ as $k$ tends to infinity, uniformly on $[\tau_1,\tau_3]$.

If
$\norm{S_2}\le p<1$, then
\begin{equation*}
\vabs{\call_2 (S_2)}(\tau) \le \frac{E(1-\ve)^2}{(1-p)^3}
\int_{\tau_1}^ \tau \frac{K_2(\tau-\sigma)\, d\sigma}{R_2^3(\sigma)},
\end{equation*}
and in consequence,
\begin{equation*}
\norm{\call_2(S_2)} \le \frac{E(1-\ve)^2 }{(1-p)^3}\sup_{\tau_1\le
\tau\le \tau_3}I(\tau).
\end{equation*}

Let us estimate $\norm{\call_2(S_2)-\call_2(\hat S_2)}$ when
$\max\bigl(\norm{S_2},\norm{\hat S_2}\bigr)$ is at most equal to $p<1$. We can
see that
\begin{equation*}
\norm{\call_2(S_2)-\call_2(\hat S_2)}\le
\frac{3E(1-\ve)^2}{(1-p)^4}\sup_{\tau_1\le 
\tau\le \tau_3}I(\tau) \norm{S_2-\hat S_2}.
\end{equation*}
Therefore, for $k$ large enough, $\call_2$ is a strict contraction
from the ball of radius $p$ about $0$ to itself.
\end{proof}

Let us prove now the estimate announced on $I$:

\begin{lemma}\label{thr:3}
The following estimate holds for $I(\tau)$ on the interval
$[\tau_1,\tau_3]$: 
\begin{equation}
I(\tau)=O\bigl(\eta^{4-3\gamma_1}+\eta^{4(1+\zeta\xi_1)}\bigr).\label{eq:24}
\end{equation}
\end{lemma}

\begin{proof}
The integral $I$ is analogous to the one defined in \eqref{cle}.

We define $\tau_2$ by
\begin{equation}
\label{taubar}
\tau_2=\frac{2\ln (\xi_2/\xi_1)}{\xi_1-\xi_2},
\end{equation}
and we consider three cases:
\begin{itemize}
\item $\tau_1\le \tau\le 2\tau_1$: in this case $H_2$ cannot be
neglected relatively to $K_2$.

\item $2\tau_1\le \tau\le \tau_2$: in this interval, the dominant term
in $\dot R_2$ will be $\dot R(\tau_1) \dot K_2(\tau-\tau_1)$ and an
elementary computation shows that this expression vanishes for
$\tau= \tau_1 +\bigl(\ln (\xi_2/\xi_1)\bigr)/\bigl(\xi_1
-\xi_2\bigr)$. Thus $\dot R_2$ crosses $0$ approximately at a time
$\tau_2/2$. 

\item $\tau_2 \le \tau\le \tau_3$: the last leg of the journey, since
$K_2$ is dominant and in $K_2$, the term involving $\exp(\xi_1(\tau
-\tau_1))$ is dominant.
\end{itemize}

Before proving these estimates, we observe that there exist
positive numbers $M$ and $m$ such that
\begin{align}
&\forall \tau\in \Er^+, \quad K_2(\tau)\le M\tau,\label{K1M}\\
&\forall\tau\in[0,\tau_2], \quad K_2(\tau)\ge m\tau.\label{K1m}
\end{align}

We tackle now the three separate sub-cases in detail.

\subsection{First interval: $\tau \in [\tau_1, 2\tau_1]$}

We remark that
\begin{equation}
R_2(\tau)\ge  R(\tau_1) H_2(\tau-\tau_1).\label{eq:16}
\end{equation}
Therefore, we can estimate $I(\tau)$ as follows:
\begin{equation}
0\le I(\tau) \le \frac{1}{ R_2(\tau)R^3(\tau_1)}\int_{\tau_1}^\tau
\frac{K_2(\tau-\sigma)}{H_2^3(\sigma-\tau_1)} \, d\sigma.\label{i11}
\end{equation}
We observe  that over $[\tau_1, 2\tau_1]$, 
\begin{equation}
H_2(\tau -\tau_1)-1=o(\tau_1),\label{eq:17}
\end{equation}
 and we use \eqref{K1M}. These observations
imply the following inequalities:
\begin{align*}
&\frac{1}{ R_2(\tau)R^3(\tau_1)}\int_{\tau_1}^\tau
\frac{K_2(\tau-\sigma)}{H_2^3(\sigma-\tau_1)} \, d\sigma\\
&\le
C\frac{1}{ R^4(\tau_1)}\int_{\tau_1}^\tau (\tau - \sigma)\, d\sigma\\
&\le C \frac{\tau_1^2}{R^4(\tau_1)}
\end{align*}
and thanks to \eqref{equiv*}, the definition 
\eqref{tauinf} of $\tau_1$ and condition \eqref{gammaetoile}, we obtain
\begin{equation}
I(\tau)=O\bigl(\eta^{4-2\gamma_1}\bigr).\label{item1}
\end{equation}

\subsection{Second interval: $\tau \in [2\tau_1, \tau_2]$}

We cut the integral $I$ into two pieces: one piece
from $\tau_1$ to $2\tau_1$ on which we work essentially as in the
previous sub-case, and a piece from $2\tau_1$ to $\tau_2$ on which we
work differently. More precisely, on $[\tau_1, 2\tau_1]$, we observe
that $R_2(\tau)\ge R(\tau_1)H_2(\tau-\tau_1)$, and on $[2\tau_1,
\tau_2]$, $R_2(\tau)\ge \dot
R(\tau_1)K_2(\tau-\tau_1)$. Therefore,
\begin{equation*}
\begin{split}
&\frac{1}{R_2(\tau)}\int_{\tau_1}^{2\tau_1}\frac{K_2(\tau-
\sigma)}{R_2^3(\sigma)}\, d\sigma\\
&\le \frac{C}{R_2(\tau)R^3(\tau_1)}\int_{\tau_1}^{2\tau_1} (\tau -
\sigma)\, d\sigma\\
&\le \frac{C \eta^{3-3\gamma_1}\tau_1\tau}{R_2(\tau)}.
\end{split}
\end{equation*}
We estimate $K_2(\tau-\tau_1)$ from below by arguing that $K_2$
increases from $0$ to a maximum, and then decreases exponentially
fast to $0$. Therefore, for all small enough $\eta$, there exists
$\tau'_1$ tending to infinity such that
$K_2(\tau'_1)=K_2(\tau_1)$. Moreover on the interval $[\tau_1,\tau'_1]$
$K_2(\tau)$ is greater than or equal to $K_2(\tau_1)$. Thus, for
all large enough $k$, $K_2(\tau)\ge K_2(\tau_1)$ on the
interval $[\tau_1,\tau'_1-\tau_1]$, and therefore
\begin{equation*}
\forall \tau\in [2\tau_1, \tau'_1],\quad R_2(\tau)\ge
\dot R(\tau_1)K_2(\tau-\tau_1)\ge\dot R(\tau_1) K_2(\tau_1).
\end{equation*}
Thus, we obtain thanks to \eqref{equiv*}
\begin{equation}
\frac{1}{R_2(\tau)}\int_{\tau_1}^{2\tau_1}\frac{K_2(\tau-
\sigma)}{R_2^3(\sigma)}\, d\sigma \le C \eta^{4-3\gamma_1}.\label{11}
\end{equation}
For the other piece, we estimate $R_2(\tau)$ from below by $\dot
R(\tau_1)K_2(\tau-\tau_1)$, and we obtain
\begin{equation}
\frac{1}{R_2(\tau)}\int_{2\tau_1}^{\tau}\frac{K_2(\tau-
\sigma)}{R_2^3(\sigma)}\,d\sigma
\le \frac{1}{R_2(\tau)\dot R^3(\tau_1)}\int_{2\tau_1}^{\tau} \frac{K_2(\tau-
\sigma)}{K_2^3(\sigma-\tau_1)}\, d\sigma.\label{ineg23}
\end{equation}
We use \eqref{K1M} and \eqref{K1m} in the integral term of the right
hand side of inequality \eqref{ineg23}, and we infer that
\begin{equation}
\begin{split}
&\frac{1}{R_2(\tau)\dot R^3(\tau_1)}\int_{2\tau_1}^{\tau} \frac{K_2(\tau-
\sigma)}{K_2^3(\sigma-\tau_1)}\, d\sigma\\
&\le \frac{C}{R_2(\tau)\dot R^3(\tau_1)} \int_{2\tau_1}^{\tau}
\frac{\tau -\sigma}{(\sigma -\tau_1)^3}\, d\sigma\\
&\le \frac{C(\tau - \tau_1)}{R_2(\tau)\dot R^3(\tau_1)\tau_1^2}\\
&\le 0(\eta^{4-2\gamma_1}).
\end{split}\label{12}
\end{equation}
The combination of \eqref{11} and \eqref{12} yields
\begin{equation}
I(\tau)=O\bigl(\eta^{4-3\gamma_1}\bigr).\label{item2}
\end{equation}

\subsection{Third interval: $\tau \in [\tau_2, \tau_3]$}

We cut now $I$ into three pieces, relative to the
subintervals $[\tau_1, 2\tau_1]$, $[2\tau_1,\tau_2]$ and
$[\tau_2,\tau]$. 

On the last two pieces, we observe that $\tau$ is far from $\tau_1$,
and we use the estimate from below
\begin{equation}
R_2(\tau)\ge \dot R(\tau_1) K_2(\tau - \tau_1).\label{eq:18}
\end{equation}
Moreover there exists $C$ such that for $\tau\ge \tau_2$ and $k$ large
enough
\begin{equation}
K_2(\tau-\tau_1) \ge C \exp\bigl(\xi_1(\tau-\tau_1)\bigr).\label{eq:19}
\end{equation}

On the first subinterval, i.e. $\sigma\in [\tau_1, 2\tau_1]$, we
use inequality \eqref{K1M}; relations \eqref{eq:16} and \eqref{eq:17}
imply that
\begin{equation*}
R_2^3(\sigma) \ge \bigl(1-o(\tau_1)\bigr) R(\tau_1)^3;
\end{equation*}
thanks to \eqref{eq:18} we can see that
\begin{equation*}
\frac{1}{R_2(\tau)} \int_{\tau_1}^{2\tau_1}
\frac{K_2(\tau-\sigma)}{R_2^3(\sigma)} \, d\sigma
\le \frac{1+o(\tau_1)}{R(\tau_1)^3 \dot R(\tau_1) K_2(\tau-\tau_1)}
\int_{\tau_1}^{2\tau_1} M(\tau-\sigma)\, d\sigma.
\end{equation*}
Thanks to \eqref{eq:19} and the asymptotics \eqref{equiv*}, we obtain
\begin{equation*}
\frac{1}{R_2(\tau)} \int_{\tau_1}^{2\tau_1}
\frac{K_2(\tau-\sigma)}{R_2^3(\sigma)} \, d\sigma \le C
\eta^{4-3\gamma_1} \exp\bigl(-\xi_1(\tau-\tau_1)\bigr) \tau_1 \tau.
\end{equation*}
Since $\tau\le \tau_3$ and $\exp(-\xi_1\tau_3)=\eta^{\zeta \xi_1}$, we
get finally
\begin{equation*}
\frac{1}{R_2(\tau)} \int_{\tau_1}^{2\tau_1}
\frac{K_2(\tau-\sigma)}{R_2^3(\sigma)} \, d\sigma \le C
\eta^{4-2\gamma_1 +\zeta \xi_1} \ln(1/\eta).
\end{equation*}

Relations \eqref{Kun} and \eqref{gammaetoile} imply that 
\begin{equation*}
4- 2\gamma_1+\zeta\xi_1>0.
\end{equation*}
We observe that 
\begin{equation}
\forall \sigma\ge 2\tau_1, \quad R_2(\sigma)\ge \dot R(\tau_1)
K_2(\sigma -\tau_1), \label{eq:25}
\end{equation}
and we use estimates \eqref{K1M} and \eqref{K1m}; therefore
\begin{equation*}
\int_{2\tau_1}^{\tau_2} \frac{K_1(\tau-\sigma)}{R_2^3(\sigma)}\,
d\sigma
\le \frac{M \tau_3}{\dot R^3(\tau_1)}\int_{2\tau_1}^{\tau_2}
\frac{d\sigma}{m^3 (\sigma -\tau_1)^3} \le \frac{M \tau_3}{2\dot
R^3(\tau_1)m^3 \tau_1^2}.
\end{equation*}
Now, thanks to \eqref{eq:18} and \eqref{eq:19}, we obtain
\begin{equation*}
\frac{1}{R_2(\tau)} 
\int_{2\tau_1}^{\tau_2}
\frac{K_1(\tau-\sigma)}{R_2^3(\sigma) } 
\, d\sigma \le
\frac{C\tau_3}{\dot R(\tau_1)^4 \tau_1^2}\exp(-\xi_1 \tau_3)
\le C \eta^{4-2\gamma_1 +\zeta\xi_1} \ln(1/\eta).
\end{equation*}

Let us consider the third piece: we use now estimate \eqref{eq:19}
on the denominator of integrand; since $K_2(\tau)\le
C\exp(\xi_1\tau)$, and thanks to \eqref{eq:25}, we have
\begin{equation*}
\begin{split}
&\frac{1}{R_2(\tau)}\int_{\tau_2}^{\tau}\frac{K_2(\tau-
\sigma)}{R_2^3(\sigma)}\,d\sigma\\
&\quad\le \frac{C}{R_2(\tau)\dot
R^3(\tau_1)}\int_{
\tau_2}^{\tau}\frac{\exp\bigl(\xi_1(\tau-\sigma)\bigr)}{\exp\bigl(3\xi_1(\sigma- 
\tau_1)\bigr)}\,d\sigma\\
&\quad\le C \eta^{4(1+\zeta\xi_1)},
\end{split}
\end{equation*}
and we conclude that the following estimate holds:
\begin{equation}\label{item3}
I(\tau)=
O\bigl(\eta^{4+\zeta\xi_1-2\gamma_1}\ln(1/\eta)+\eta^{4(1+\zeta\xi_1)}
\bigr). 
\end{equation}
We have to keep the two terms in the above expression, since we have
no way to ascertain the order of the exponents of $\eta$.

When we compare the exponents in \eqref{item1} and \eqref{item2}, we
find that the exponent of $\eta$ in \eqref{item2} is the smaller; when
we look at the exponents in \eqref{item3} to the exponent in \eqref{item2}
we find that $4+\zeta \xi_1 -2\gamma_1$ is strictly larger than
$4-3\gamma_1$, and this leads to the conclusion \eqref{eq:24}.
\end{proof}

We state now the main result of this section:

\begin{proposition}
The following estimates hold:
\begin{align}
&R(\tau)\sim R_2(\tau)\text{ uniformly over $[\tau_1,\tau_3]$},
\label{equivR1}\\
&R(\tau_3) \sim \frac{\dot s(0)}{2\sqrt\Delta}\eta^{-(1+\zeta
\xi_1)},\label{eq:27}\\
&\dot R(\tau_3) \sim \frac{\dot s(0)}{2\sqrt\Delta} \xi_1 
\eta^{-(1+\zeta \xi_1)}.\label{eq:28}
\end{align}
\end{proposition}

\begin{proof}
Theorem \ref{contractionbis} implies the uniform
equivalence~\eqref{equivR1}, and \eqref{eq:27} is an immediate
consequence of \eqref{equivR1}.

Let us prove an estimate of the derivative $\dot R$ at $\tau_3$:
\begin{equation}
\dot R(\tau_3)\sim \dot R_2(\tau_3)\sim \frac{\dot
s(0)}{2\sqrt\Delta}\xi_1\eta^{-(1+\zeta\xi_1)}.\label{equivdotR1} 
\end{equation}
We observe that
\begin{equation*}
\dot R(\tau_3)=\dot R_2(\tau_3)+\int_{\tau_1}^{\tau_3}\frac{\partial
K_2}{\partial \tau}(\tau_3
-\sigma)\frac{E(1-\ve)^2}{(R_2+S_2)(\sigma)^3}\,d\sigma.
\end{equation*}
Therefore,
\begin{equation*}
\abs{\dot R(\tau_3) - \dot R_2(\tau_3)}\le C \int_{\tau_1}^{\tau_3}
\vabs{\frac{\partial 
K_2}{\partial \tau}(\tau_3 -\sigma)}\frac{1}{R_2^3(\sigma)}\, d\sigma.
\end{equation*}
There exists a constant $C$ such that for all $\sigma\ge 0$
\begin{equation*}
\vabs{\frac{\partial K_2}{\partial \tau}(\sigma)}\le C\exp(\xi_1\sigma).
\end{equation*}
We use the method which gave estimate \eqref{item3}: we cut the
integration interval into the three subintervals $[\tau_1, 2\tau_1]$,
$[2\tau_1,\tau_2]$ and $[\tau_2,\tau_3]$, and on each of these
subintervals, we estimate $R_2$ from below exactly as in this
calculation. Details are left to the reader, and we obtain
\begin{equation*}
\int_{\tau_1}^{\tau_3} \frac{\exp(\xi_1\sigma)}{R_2^3(\sigma)}\,
d\sigma = O\bigl(\eta^{3-2\gamma_1-\xi_1\zeta} +
\eta^{3(1+\xi_1\zeta)}\bigr). 
\end{equation*}
The equivalent of $\dot R_2(\tau_3)$ is obtained immediately from the
explicit formula \eqref{R1} for $R_2$ and the equivalents \eqref{equiv*}.
Hence we infer that \eqref{equivdotR1} holds.
\end{proof}

\section{The final asymptotics}\label{sec:final-asymptotics}

In this section, we show that for large enough times $R(\tau)$ is
bounded from above. In view of \eqref{coin:eq:10}, this estimate will
enable us to show that the angular velocity is bounded from below, and
hence, the polar angle $\Theta$ will cross through $\bar \theta$.

\begin{theorem}
There exists a strictly
positive number $R_M$ and a time $\tau_4$ such that 
\begin{equation*}
\forall \tau\ge \tau_4,\quad R(\tau)\le R_M.
\end{equation*}
\end{theorem}

\begin{proof}
Denote
\begin{equation*}
x=\begin{pmatrix}
R\\ \dot R
\end{pmatrix}, \quad N(x)=\begin{pmatrix}
0\\ E(1-\ve)^2/R^3
\end{pmatrix}, \quad M=\begin{pmatrix}
0&1\\ -1&-2\alpha
\end{pmatrix}.
\end{equation*}

With these notations, equation \eqref{eqrnsc2} can be written
\begin{equation}
\dot x=Mx + N(x).\label{systeme}
\end{equation}
We observe that in the domain $(0,\infty)\times \Er$, \eqref{systeme}
has exactly one critical point at
\begin{equation*}
x_c =\begin{pmatrix}
R_c\\0
\end{pmatrix}, \quad R_c=\left(E(1-\ve)^2\right)^{1/4}.
\end{equation*}
This critical point is attractive, as an examination of the
linearization of \eqref{systeme} around $x_c$ shows. Moreover, there
is a Lyapunov functional given by
\begin{equation}
F(x)=x_1^2 +\frac{E(1-\ve)^2}{x_1^2} + x_2^2.\label{lyapunov}
\end{equation}
Therefore, given $x(\tau)$ with $x_1(\tau)>0$, we can see that for all
$\tau' \ge \tau$, 
$F(x(\tau'))$ is at most equal to $F(x(\tau))$, and in particular,
$x(\tau')$ remains bounded. We see that when $\tau$ tends to
infinity, $x(\tau)$ tends to the critical point $x_c$.

The spectrum of $M$ is $\{\xi_1,\xi_2\}$; therefore, the matrix
\begin{equation*}
Q=\int_0^\infty \exp(sM^*)\exp(sM)\, ds
\end{equation*}
is well defined, symmetric, positive and definite. In particular,
if $\lambda_1$ is the smallest eigenvalue of $Q$ and $\lambda_2$ is the
largest eigenvalue of $Q$,
\begin{equation}
\lambda_1 x^* x\le x^*Q x\le \lambda_2 x^* x.\label{eq:29}
\end{equation}
If we let $x(t)=e^{tM}x_0$, we observe that
\begin{equation*}
\begin{split}
\frac{d}{dt} x(t)^*Q x(t) &=\frac{d}{dt} \int_0^\infty
x_0^*\exp(tM^*)\exp(sM^*)\exp(sM) \exp(tM)x_0\, ds\\
&=\frac{d}{dt} \int_t^\infty x_0^*\exp(sM^*)\exp(sM)x_0\, ds\\
&=-x(t)^*x(t) \le -\lambda_2^{-1} x(t)^*Q x(t).
\end{split}
\end{equation*}
As $x_0\in \Er^2$ and $t\ge 0$ are arbitrary in the above calculation,
we have proved indeed that for all $x\in \Er^2$
\begin{equation}
2x^*QMx \le -\lambda_2^{-1} x^*Qx.\label{quad}
\end{equation}

Since $\dot x=Mx +N(x)$ we have the inequality
\begin{equation*}
\frac{d}{d\tau} x^* Qx \le -\frac{x^* Qx}{\lambda_2} +2\bigl(x^* Q
x\bigr)^{1/2} \sqrt{\lambda_2} \vabs{N(x)}.
\end{equation*}
We seek a number $\bar R$ such that if $R(\tau)\ge \bar R$, then
\begin{equation}
\frac{d}{d\tau} \bigl(x^* Qx\bigr) \le -\frac{x^*
Qx}{2\lambda_2}.\label{eq:26} 
\end{equation}
Indeed, in order to satisfy \eqref{eq:26}, it suffices to have
\begin{equation*}
2\sqrt{\lambda_2} \abs{N(x)}\bigl(x^* Q x\bigr)^{1/2}\le \frac{x^*
Qx}{2\lambda_2}, 
\end{equation*}
or equivalently,
\begin{equation*}
\abs{N(x)} \le \frac{(x^* Qx)^{1/2}}{4 \lambda_2^{3/2}}.
\end{equation*}
But $\abs{N(x)} =E(1-\ve)^2/R^3$ and $\abs{x}\ge R$, so that, with the
help of \eqref{eq:29}, it suffices to satisfy
\begin{equation*}
\frac{E}{\bar R^3} \le \frac{\lambda_1^{1/2}\bar R}{4 \lambda_2^{3/2}},
\end{equation*}
i.e. 
\begin{equation*}
\bar R^4 \ge \frac{4\lambda_2^{3/2} E}{\lambda_1^{1/2}}.
\end{equation*}

We shall show now that if we choose $\bar R$ such that
\begin{equation*}
\bar R >
\max\left(\left(\frac{4E\lambda_2^{3/2}}{\lambda_1^{1/2}}\right), R_c\right),
\end{equation*}
then there exists $\tau_4$ such that
\begin{equation}
R(\tau_4)= \bar R.\label{eq:30}
\end{equation}
Indeed, we know from \eqref{eq:27} that $R(\tau_3)\gg 1$, and that
the limit of $R(\tau)$ as $\tau$ tends to infinity is $R_c$;
therefore, $R(\tau)$ must cross $\bar R$. We denote by $\tau_4$ the
smallest time in $[\tau_3,\infty)$ such that \eqref{eq:30} holds.

On the interval $[\tau_3,\tau_4]$, the differential inequality
\eqref{eq:26} implies
\begin{equation*}
\bigl(x^* Qx\bigr)(\tau_4)\le \bigl(x^* Qx\bigr)(\tau_3)
\exp\bigl(-(\tau_4-\tau_3)/2\lambda_2\bigr),
\end{equation*}
whence
\begin{equation*}
\frac{\tau_4 -\tau_3}{2\lambda_2} \le\ln\bigl(x^*Qx\bigr)(\tau_3)
-\ln\bigl(x^*Qx\bigr) (\tau_4).
\end{equation*}
But $\bigl(x^*Qx\bigr)(\tau_4)\ge \lambda_1 \bar R^2$, and we obtain
the inequality
\begin{equation*}
\tau_4 \le \tau_3 +2\lambda_2\bigr[\ln \bigl(x^*Qx\bigr)(\tau_3)
-\ln\bigl(\lambda_1 \bar R^2\bigr)\bigr].
\end{equation*}
In particular, there exists $C$ such that
\begin{equation*}
\tau_4\le \tau_3 + C\ln(1/\eta).
\end{equation*}

We also need an estimate on $\dot R(\tau_4)$. We first show that it is
less than or equal to $0$. By \eqref{eq:28} we know that
$\dot R(\tau_3)<0$. Denote by $(\tau_3,\tau_5)$ the connected
component of $\{\tau >\tau_3:\dot R(\tau)<0\}$ whose boundary contains
$\tau_3$. 

If $\tau_5=\infty$, it is clear that $\dot R(\tau_4)\le 0$. Assume
that $\tau_5<\infty$ and that $\dot R(\tau_4)>0$; then $\tau_5<\tau_4$
and $\dot R(\tau_5)$ vanishes.

We infer from differential equation~\eqref{eqrnsc2} that
\begin{equation*}
\ddot R(\tau_5)=-R(\tau_5) +\frac{E(1-\ve)^2}{R(\tau_5)^3},
\end{equation*}
but $R(\tau_5)>R_c$, because $R(\tau_5)\ge R(\tau_4)= \bar R$;
therefore
\begin{equation}
\ddot R(\tau_5)<0.\label{eq:31}
\end{equation}
On the other hand, as $\dot R(\tau)$ is negative on $(\tau_3,\tau_5)$
and vanishes at $\tau_5$, a straightforward sign argument shows that
\begin{equation*}
\ddot R(\tau_5) \ge 0,
\end{equation*}
which contradicts \eqref{eq:31}.

Now, we prove that
\begin{equation*}
\dot R(\tau_4)\ge \xi_1 R(\tau_4).
\end{equation*}
This will be a consequence of the following inequality
for all
$\tau\ge\tau_1$ and for all large enough $k$:
\begin{equation}
\label{plus}
\dot R(\tau)-\xi_1 R(\tau)\ge 0.
\end{equation}
We observe that
\begin{equation*}
\frac{d}{d\tau}\bigl( \dot R-\xi_1 R\bigr) =\xi_2\bigl(R-\xi_1 R\bigr)
+\frac{E(1-\ve)^2}{R^3}. 
\end{equation*}
When we integrate this differential relation, we find that
\begin{equation*}
\begin{split}
\bigl( \dot R-\xi_1 R\bigr)(\tau)&=\exp\bigl(\xi_2(\tau-\tau_1)\bigr)
\bigl( \dot R-\xi_1 R\bigr)(\tau_1)\\
&\quad +\int_{\tau_1}^\tau
\exp\bigl(\xi_2(\tau-\sigma) \bigr)\frac{E(1-\ve)^2}{R^3(\sigma)}\,
d\sigma.
\end{split}
\end{equation*}
For $k$ large enough, the equivalences \eqref{equiv*} show that
$\bigl( \dot R-\xi_1 R\bigr)(\tau_1)$ is strictly positive, and
\eqref{plus} follows immediately.

We infer now from \eqref{plus} and the sign condition on $\dot
R(\tau_4)$ that
\begin{equation*}
F(x(\tau_4))\le \bar F= \bar R^2 + E(1-\ve)^2\bar R^{-2} + \xi_1^2\bar
R^2. 
\end{equation*}
Since the Lyapunov functional decreases along trajectories of the
system, we obtain for all $\tau\ge \tau_4$ the inequalities
\begin{equation}
\frac{E(1-\ve)}{\bar F}\le R(\tau)^2 \le \bar F, \quad \dot R(\tau)^2
\le \bar F.\label{eq:32}
\end{equation}
This is the final estimate we needed before the conclusion.
\end{proof}

We can now state the following corollary relative to the existence of
the time $\bar \tau$:

\begin{corollary}
There exists a time $\bar \tau\in(0,\infty)$ such that $\Theta(\bar
\tau)=\bar \theta$.
\end{corollary}

\begin{proof}
We know from \eqref{coin:eq:10} that $\Theta$ is an increasing
function of $\tau$; if there is a time $\bar \tau\le \tau_4$ for which
$\Theta(\bar \tau)=\bar \theta$, the conclusion is clear. Assume
otherwise; then, with the notations of \eqref{eq:32}, we can see that
\begin{equation*}
\dot \Theta(\tau) =\frac{\sqrt E(1-\ve)}{R^2} \ge \frac{\sqrt
E(1-\ve)}{\bar F},
\end{equation*}
and the conclusion is also clear.
\end{proof}

\section{The case $\bar \theta \ge \pi/2$}

In this section we estimate from below the first time $\bar \tau$ at
which $\Theta(\bar \tau)=\bar \theta$; we expect that $\bar \tau$ will
be comparable to $\tau_1$, but this is not correct. Recall the
definition \eqref{valeps} of $\ve$; in this definition, the exponent
of $\eta$ is 
\begin{equation*}
4\sqrt{\Delta}/\abs{\xi_1};
\end{equation*}
define a number $r$ by
\begin{equation}
r=\min\bigl(\gamma_1, 4\sqrt{\Delta}/\vabs{\xi_1}\bigr).\label{coin:eq:11}
\end{equation}

Now, we can state the following theorem:

\begin{theorem}
If $\bar \theta >\pi/2$, then for large enough $k$, $\bar  \tau\ge
\tau_1$; if $\bar \theta =\pi/2$, then there exists a strictly 
positive number $C$ such that
\begin{equation}
\bar \tau \ge C \eta^{\max(2-r, \gamma_1)}.\label{coin:eq:4}
\end{equation}
\end{theorem}

\begin{proof}
We argue as follows: assume $\bar \tau\le\tau_1$; we recall
estimate~\eqref{cadre}:
\begin{equation*}
\forall \tau\in[0,\tau_1], \quad 
(1-\beta)\Theta_1(\tau) \le \Theta(\tau) \le (1+\beta) \Theta_1(\tau),
\end{equation*}
where $\beta$ is given by \eqref{coin:eq:1} and $\norm{S_1}=p$
satisfies \eqref{eq:23}. The assumption \eqref{gammaetoile} implies
$4\gamma_1 -4<\gamma_1$, and thus \eqref{eq:23} simplifies as
\begin{equation*}
p=O\bigl(\eta^{\gamma_1} +\ve\bigr).
\end{equation*}
The definition \eqref{coin:eq:11} of $r$ implies that
\begin{equation*}
\beta=O\bigl(\eta^r\bigr).
\end{equation*}
Moreover, relation \eqref{formTheta} leads to 
\begin{equation*}
\Theta_1(\tau)=\frac{\pi}{2} -\arctan\frac{\sqrt{E}}{W(\tau-\tau_0)}
+\arctan \frac{W\tau_0}{\sqrt{E}}.
\end{equation*}
This relation implies immediately that
\begin{equation*}
\lim_{k\to\infty} \Theta_1(\tau_1)=\frac{\pi}{2},
\end{equation*}
and therefore, thanks to \eqref{cadre}
\begin{equation*}
\lim_{k\to\infty} \Theta(\tau_1)=\frac{\pi}{2}.
\end{equation*}
If $\bar \theta >\pi/2$, the last relation implies immediately that
for $k$ large enough, $\bar \tau$ is at least equal to $\tau_1$. 

Assume now that $\bar \theta=\pi/2$; now, the situation is more
delicate, since none of the inequalities established so far implies an
estimate on $\bar \tau$. If $\bar \tau\ge \tau_1$, we are done.
Otherwise, we shall estimate $\bar \tau$ from below. Already, relation
\eqref{cadre} implies
\begin{equation*}
\Theta_1(\bar \tau) \ge\frac{\pi}{2(1+C\eta^r)},
\end{equation*}
or in other words
\begin{equation*}
\frac{\pi}{2} -\arctan \frac{\sqrt E}{W(\bar \tau -\tau_0)} +\arctan
\frac{W\tau_0}{\sqrt E} \ge\frac{\pi}{2(1+ C\eta^r)},
\end{equation*}
which implies
\begin{equation*}
\arctan \frac{\sqrt E}{W(\bar \tau-\tau_0)} \le O\bigl(\eta^r\bigr),
\end{equation*}
and therefore
\begin{equation*}
\bar \tau-\tau_0 \ge C\eta^{2-r}.
\end{equation*}
Thus, we have shown that
\begin{equation*}
\bar \tau \ge C \eta^{2-r}.
\end{equation*}
If $\gamma_1> 2-r$, the relations
\begin{equation*}
\tau_1 =\eta^{\gamma_1} \ge \bar \tau \ge C \eta^{2-r}
\end{equation*}
are contradictory for $k$ large; therefore
\begin{equation*}
\gamma_1 > 2-r, \text{ $k$ large} \Longrightarrow \bar \tau \ge \tau_1.
\end{equation*}
Thus, we have shown
~\eqref{coin:eq:4}.
\end{proof}

We deduce the following estimates from \eqref{coin:eq:4} and the
asymptotics of sections \ref{sec:earliest},
\ref{sec:second}~ and \ref{sec:final-asymptotics} 
\begin{equation}
\begin{split}
&C\eta^{\max(1-r,\gamma_1-1)} \le R(\bar \tau) \le
\frac{C'}{\eta},\\ 
&\abs{\dot R(\bar \tau)}\le\frac{C}{\eta}
\end{split}\label{coin:eq:5}
\end{equation}

We are able to show now the main result of this section:

\begin{theorem}\label{coin:thr:1}
For $\bar \theta \ge \pi/2$, as $k$ tends to infinity, $u_k$ converges
uniformly on the compact sets of $\Er^+$ to $u_\infty$ given by
\begin{equation*}
u_\infty(t)=\begin{cases}
u(0) + t\Pi_1 \dot u(0) &\text{if $0 \le t \le t_0$,}\\
0 &\text{if $t_0 \le t$,}
\end{cases}
\end{equation*}
where $\Pi_1$ is the projection on the line $\{x_1=0\}$; see
Fig. \ref{fig:region}. 
\end{theorem}

\begin{proof}
We go back to the original scales and time $\bar t= t_0 + \bar
\tau/\sqrt k$; then
\begin{equation*}
u_k(\bar t)=r(\bar t) e^{i\bar \theta}, \quad \dot u_k(\bar t) =
\bigl(\dot r(\bar t + ir(\bar t) \dot \theta(\bar t)\bigr) 
e^{i\bar\theta}. 
\end{equation*}
Therefore, in coordinates $y_1, y_2$ (see Fig. \ref{fig:region}), we
have the relations
\begin{equation}
\begin{split}
y_1(\bar t)&= \eta R(\bar \tau)/\sqrt k,\\
y_2(\bar t)&=0,\\
\dot y_1(\bar t) &= \eta \dot R(\bar \tau),\\
\dot y_2(\bar t)&= \eta (1-\ve) \sqrt{E}\, /R(\bar \tau).
\end{split}\label{coin:eq:7}
\end{equation}
We have also the estimate
\begin{equation}
\dot R(\bar \tau) -\xi_1 R(\bar \tau)\ge 0.\label{coin:eq:8}
\end{equation}
If $\bar \tau \ge \tau_1$, \eqref{coin:eq:8} is a consequence of
\eqref{plus}. Otherwise, we observe that $\bar \tau$ belongs to
$[\eta^3,\tau_1]$ for all large enough $k$; therefore, we are able to
use the equivalences \eqref{eq:21} and \eqref{eq:22}, whence
\begin{equation*}
\dot R(\bar \tau) -\xi_1 R(\bar \tau) \sim \dot R_1(\bar \tau) -\xi_1
R_1(\bar \tau),
\end{equation*}
which is valid because the dominant term in the right hand side of the
above expression does not vanish; indeed, the expression \eqref{solR}
of $R_1$ and \eqref{solR'} of $\dot R_1$, we can see that
\begin{equation*}
\dot R_1(\bar \tau)-\xi_1 R_1(\bar \tau) \sim \sqrt{W} -\xi_1
\sqrt{W}\, \bar \tau \sim C \eta^{-1},
\end{equation*}
which implies \eqref{coin:eq:8}; in the original coordinates,
\eqref{coin:eq:8} translates as
\begin{equation}
\dot y_1(\bar t) -\sqrt k\, \xi_1 y_1(\bar t)\ge 0.\label{coin:eq:9}
\end{equation}
We infer from estimate \eqref{coin:eq:5} that
\begin{equation*}
\begin{split}
y_1(\bar t) &=O(1/\sqrt k),\\
\dot y_1(\bar t)&=O(1),\\
\dot y_2(\bar t) &=O\bigl(\eta^{1-\max(1-r, \gamma_1-1)}\bigr) =o(1).
\end{split}
\end{equation*}
In the coordinates $y_1$ and $y_2$, the system \eqref{approx} can be
rewritten 
\begin{align}
&\ddot y_1 + 2\alpha \sqrt {k}\, \dot y_1 + k y_1 =0 \label{coin:eq:12}\\
\intertext{as long as $y_1\ge 0$ and}
&\ddot y_2 =0.\notag
\end{align}
But the explicit solution of \eqref{coin:eq:12} with initial data
~\eqref{coin:eq:7} is given by
\begin{equation*}
\begin{split}
y_1(t)&=\dot y_1(\bar t) \frac{\exp\bigl(\xi_1(t-\bar t)\sqrt k\bigr)
-\exp(\bigl( \xi_2(t-\bar t)\sqrt k\bigr)}{2\sqrt{\Delta k}} \\
&\quad + y_1(\bar t) \frac{
\xi_1 \exp\bigl(\xi_2(t-\bar t)\sqrt k\bigr)- \xi_2
\exp\bigl(\xi_1(t-\bar t)\sqrt k\bigr)}{2\sqrt \Delta}.
\end{split}
\end{equation*}
If $\dot y(t_1)$ is non negative, it is clear that $y_1$ stays non
negative for all time larger than $\bar t$.
If $\dot y_1(\bar t)$ is negative, we use \eqref{coin:eq:9}: we
estimate from below $\dot y_1(\bar t)$ by $\sqrt k \, \xi_1 y_1(\bar
t)$, and after simplifications, we get
\begin{equation*}
y_1(t)\ge y_1(\bar t) \frac{\xi_1 -\xi_2}{2\sqrt\Delta}
\exp\bigl(\xi_1(t-\bar t) \sqrt k\bigr).
\end{equation*}
Therefore, \eqref{coin:eq:12} holds for all $t\ge \bar t$. In
particular, 
\begin{equation*}
\forall t\ge \bar t,\quad\vabs{y_1(t)} =O\bigl(1/\sqrt k\bigr),
\end{equation*}
and
\begin{equation*}
 \forall
t\ge \bar t,\quad y_2(t)=O\bigl((t-\bar t) \eta^{1-\max(1-r,
\gamma_1-1)}\bigr),
\end{equation*}
which proves that in this case the limit of $y_1$ and $y_2$ is $0$,
as $k$ tends to infinity.
\end{proof}

\bibliographystyle{plain}
\bibliography{../impact}

\end{document}